# On the false discovery rates of a frequentist: Asymptotic expansions

## Anirban DasGupta[1] Tonglin Zhang[2]

*Purdue University*

**Abstract:** Consider a testing problem for the null hypothesis $H_0 : \theta \in \Theta_0$. The standard frequentist practice is to reject the null hypothesis when the p-value is smaller than a threshold value $\alpha$, usually 0.05. We ask the question how many of the null hypotheses a frequentist rejects are actually true. Precisely, we look at the Bayesian false discovery rate $\delta_n = P_g(\theta \in \Theta_0 | p - value < \alpha)$ under a proper prior density $g(\theta)$. This depends on the prior $g$, the sample size $n$, the threshold value $\alpha$ as well as the choice of the test statistic. We show that the Benjamini–Hochberg FDR in fact converges to $\delta_n$ almost surely under $g$ for any fixed $n$. For one-sided null hypotheses, we derive a third order asymptotic expansion for $\delta_n$ in the continuous exponential family when the test statistic is the MLE and in the location family when the test statistic is the sample median. We also briefly mention the expansion in the uniform family when the test statistic is the MLE. The expansions are derived by putting together Edgeworth expansions for the CDF, Cornish–Fisher expansions for the quantile function and various Taylor expansions. Numerical results show that the expansions are very accurate even for a small value of $n$ (e.g., $n = 10$). We make many useful conclusions from these expansions, and specifically that the frequentist is not prone to false discoveries except when the prior $g$ is too spiky. The results are illustrated by many examples.

## 1. Introduction

In a strikingly interesting short note, Sorić [19] raised the question of establishing upper bounds on the proportion of fictitious statistical discoveries in a battery of independent experiments. Thus, if $m$ null hypotheses are tested independently, of which $m_0$ happen to be true, but $V$ among these $m_0$ are rejected at a significance level $\alpha$, and another $S$ among the false ones are also rejected, Sorić essentially suggested $E(V)/(V + S)$ as a measure of the false discovery rate in the chain of $m$ independent experiments. Benjamini and Hochberg [3] then looked at the question in much greater detail and gave a careful discussion for what a correct formulation for the false discovery rate of a group of frequentists should be, and provided a concrete procedure that actually physically controls the groupwise false discovery rate. The problem is simultaneously theoretically attractive, socially relevant, and practically important. The practical importance comes from its obvious relation to statistical discoveries made in clinical trials, and in modern microarray experiments. The continued importance of the problem is reflected in two recent articles, Efron [7], and Storey [21], who provide serious Bayesian connections and advancements in the problem. See also Storey [20], Storey, Taylor and Siegmund [22], Storey and Tibshirani [23], Genovese and Wasserman [10], and Finner and Roters [9], among many others in this currently active area.

[1]Department of Statistics, Purdue University, 150 North University Street, West Lafayette, IN 47907-2067, e-mail: `dasgupta@stat.purdue.edu`

[2]Department of Statistics, Purdue University, 150 North University Street, West Lafayette, IN 47907-2067, e-mail: `tlzhang@stat.purdue.edu`

Around the same time that Sorić raised the issue of fictitious frequentist discoveries made by a mechanical adoption of the use of p-values, a different debate was brewing in the foundation literature. Berger and Sellke [2], in a thought provoking article, gave analytical foundations to the thesis in Edwards, Lindman and Savage [6] that the frequentist practice of rejecting a sharp null at a traditional 5% level amounts to a rush to judgment against the null hypothesis. By deriving lower bounds or exact values for the minimum value of the posterior probability of a sharp null hypothesis over a variety of classes of priors, Berger and Sellke [2] argued that p-values traditionally regarded as small understate the plausibility of nulls, at least in some problems. Casella and Berger [5], gave a collection of theorems that show that the discrepancy disappears under broad conditions if the null hypothesis is composite one-sided. Since the articles of Berger and Sellke [2] and Casella and Berger [5], there has been an avalanche of activity in the foundation literature on the safety of use of p-values in testing problems. See Hall and Sellinger [12], Sellke, Bayarri and Berger [18], Marden [14] and Schervish [17] for a contemporary exposition.

It is conceptually clear that the frequentist FDR literature and the foundation literature were both talking about a similar issue: is the frequentist practice of rejecting nulls at traditional p-values an invitation to rampant false discoveries? The structural difference was that the FDR literature did not introduce a formal prior on the unknown parameters, while the foundation literature did not go into multiple testing, as is the case in microarray or other emerging interesting applications. The purpose of this article is to marry the two schools together, while giving a new rigorous analysis of the interesting question: "how many of the null hypotheses a frequentist rejects are actually true" and the flip side of that question, namely, "how many of the null hypotheses a frequentist accepts are actually falses". The calculations are completely different from what the previous researchers have done, although we then demonstrate that our formulation directly relates to *both* the traditional FDR calculations, and the foundational effort in Berger and Sellke [2], and others. We have thus a dual goal; providing a new approach, and integrating it with the two existing approaches.

In Section 2, we demonstrate the connection in very great generality, without practically any structural assumptions at all. This was comforting. As regards to concrete results, it seems appropriate to look at the one parameter exponential family, it being the first structured case one would want to investigate. In Section 3, we do so, using the MLE as the test statistic. In Section 4, we look at a general location parameter, but using the median as the test statistic. We used the median for two reasons. First, for general location parameters, the median is credible as a test statistic, while the mean obviously is not. Second, it is important to investigate the extent to which the answers depend on the choice of the test statistic; by studying the median, we get an opportunity to compare the answers for the mean and the median in the special normal case.To be specific, let us consider the one sided testing problem based on an i.i.d. sample $X_1, \ldots, X_n$ from a distribution family with parameter $\theta$ in the parameter space $\Omega$ which is an interval of $R$. Without loss of generality, we assume $\Omega = (\underline{\theta}, \bar{\theta})$ with $-\infty \leq \underline{\theta} < \bar{\theta} \leq \infty$. We consider the testing problem

$$H_0 : \theta \leq \theta_0 \text{ vs } H_1 : \theta > \theta_0,$$

where $\theta_0 \in (\underline{\theta}, \bar{\theta})$. Suppose the $\alpha$, $0 < \alpha < 1$, level test rejects $H_0$ if $T_n \in C$, where $T_n$ is a test statistic. We study the behavior of the quantities,

$$\delta_n = P(\theta \leq \theta_0 | T_n \in C) = P(H_0 | p - value < \alpha)$$

and

$$\epsilon_n = P(\theta > \theta_0 | T_n \notin C) = P(H_1 | p - value \geq \alpha).$$

Note that $\delta_n$ and $\epsilon_n$ are inherently Bayesian quantities. By an almost egregious abuse of nomenclature, we will refer to $\delta_n$ and $\epsilon_n$ as type I and type II errors in this article. Our principal objective is to obtain third order asymptotic expansions for $\delta_n$ and $\epsilon_n$ assuming a Bayesian proper prior for $\theta$. Suppose $g(\theta)$ is any sufficiently smooth proper prior density of $\theta$. In the regular case, the expansion for $\delta_n$ we obtain is like

$$(1) \qquad \delta_n = \frac{P(\theta \leq \theta_0, T_n \in C)}{P(T_n \in C)} = \frac{c_1}{\sqrt{n}} + \frac{c_2}{n} + \frac{c_3}{n^{3/2}} + O(n^{-2}),$$

and the expansion for $\epsilon_n$ is like

$$(2) \qquad \epsilon_n = \frac{P(\theta > \theta_0, T_n \notin C)}{P(T_n \notin C)} = \frac{d_1}{\sqrt{n}} + \frac{d_2}{n} + \frac{d_3}{n^{3/2}} + O(n^{-2}),$$

where the coefficients $c_1, c_2, c_3$, $d_1$, $d_2$, and $d_3$ depend on the problem, the test statistic $T_n$, the value of $\alpha$ and the prior density $g(\theta)$. In the nonregular case, the expansion differs qualitatively; for both $\delta_n$ and $\epsilon_n$ the successive terms are in powers of $1/n$ instead of the powers of $1/\sqrt{n}$. Our ability to derive a third order expansion results in a surprisingly accurate expansion, sometimes for $n$ as small as $n = 4$. The asymptotic expansions we derive are not just of theoretical interest; the expansions let us conclude interesting things, as in Sections 3.2 and 4.5, that would be impossible to conclude from the exact expressions for $\delta_n$ and $\epsilon_n$.

The expansions of $\delta_n$ and $\epsilon_n$ require the expansions of the numerators and the denominators of (1) and (2) respectively. In the regular case, the expansion of the numerator of (1) is like

$$(3) \qquad A_n = P(\theta \leq \theta_0, T_n \in C) = \frac{a_1}{\sqrt{n}} + \frac{a_2}{n} + \frac{a_3}{n^{3/2}} + O(n^{-2})$$

and the expansion of the numerator of (2) is like

$$(4) \qquad \tilde{A}_n = P(\theta > \theta_0, T_n \notin C) = \frac{\tilde{a}_1}{\sqrt{n}} + \frac{\tilde{a}_2}{n} + \frac{\tilde{a}_3}{n^{3/2}} + O(n^{-2}).$$

Then, the expansion of the denominator of (1) is

$$(5) \qquad B_n = P(T_n \in C) = A_n + \lambda - \tilde{A}_n = \lambda - \frac{b_1}{\sqrt{n}} - \frac{b_2}{n} - \frac{b_3}{n^{3/2}} + O(n^{-2}),$$

where $\lambda = P(\theta > \theta_0) = \int_{\theta_0}^{\bar{\theta}} g(\theta) d\theta$ and assume $0 < \lambda < 1$, $b_1 = \tilde{a}_1 - a_1$, $b_2 = \tilde{a}_2 - a_2$ and $b_3 = \tilde{a}_3 - a_3$, and the expansion of the denominator of (2) is

$$(6) \qquad \tilde{B}_n = P(T_n \notin C) = 1 - B_n = 1 - \lambda + \frac{b_1}{\sqrt{n}} + \frac{b_2}{n} + \frac{b_3}{n^{3/2}} + O(n^{-2}).$$

Then, we have

$$c_1 = \frac{a_1}{\lambda}, c_2 = \frac{a_1 b_1}{\lambda^2} + \frac{a_2}{\lambda}, c_3 = \frac{a_3}{\lambda} + \frac{a_1 b_2 + a_2 b_1}{\lambda^2} + \frac{a_1 b_1^2}{\lambda^3},$$

(7)

$$d_1 = \frac{\tilde{a}_1}{1-\lambda}, d_2 = \frac{\tilde{a}_2}{1-\lambda} - \frac{\tilde{a}_1 b_1}{(1-\lambda)^2}, d_3 = \frac{\tilde{a}_3}{1-\lambda} - \frac{\tilde{a}_2 b_1 + \tilde{a}_1 b_2}{(1-\lambda)^2} + \frac{\tilde{a}_1 b_1^2}{(1-\lambda)^3}.$$

We will frequently use the three notations in the expansions: the standard normal PDF $\phi$, the standard normal CDF $\Phi$ and the standard normal upper $\alpha$ quantile $z_\alpha = \Phi^{-1}(1-\alpha)$.

The principal ingredients of our calculations are Edgeworth expansions, Cornish–Fisher expansions and Taylor expansions. The derivation of the expansions became very complex. But in the end, we learn a number of interesting things. We learn that typically the false discovery rate $\delta_n$ is small, and smaller than the pre-experimental claim $\alpha$ for quite small $n$. We learn that typically $\epsilon_n > \delta_n$, so that the frequentist is less vulnerable to false discovery than to false acceptance. We learn that only priors very spiky at the boundary between $H_0$ and $H_1$ can cause large false discovery rates. We also learn that these phenomena do not really change if the test statistic is changed. So while the article is technically complex and the calculations are long, the consequences are rewarding. The analogous expansions are qualitatively different in the nonregular case. We could not report them here due to shortage of space. We should also add that we leave open the question of establishing these expansions for problems with nuisance parameters, multivariate problems, and dependent data. Results similar to ours are expected in such problems.

## 2. Connection to Benjamini and Hochberg, Storey and Efron's work

Suppose there are $m$ groups of iid samples $X_{i1}, \ldots, X_{in}$ for $i = 1, \ldots, m$. Assume $X_{i1}, \ldots, X_{in}$ are iid with a common density $f(x, \theta_i)$, where $\theta_i$ are assumed iid with a CDF $G(\theta)$ which does not need to have a density in this section. Then, the prior $G(\theta)$ connects our Bayesian false discovery rate $\delta_n$ to the usual frequentist false discovery rate. In the context of our hypothesis testing problem, the frequentist false discovery rate, which has been recently discussed by Benjamini and Hochberg [3], Efron [7] and Storey [21], is defined as

$$FDR = FDR(\theta_1, \ldots, \theta_m) = E_{\theta_1, \ldots, \theta_m} \left\{ \frac{\sum_{i=1}^m I_{T_{ni} \in C, \theta_i \le \theta_0}}{(\sum_{i=1}^m I_{T_{ni} \in C}) \vee 1} \right\},$$

(8)

where $T_{ni}$ is the test statistic based on the samples $X_{i1}, \ldots, X_{in}$. It will be shown below that for any fixed $n$ as $m \to \infty$, the frequentist false discovery rate $FDR$ goes to the Bayesian false discovery rate $\delta_n$ almost surely under the prior distribution $G(\theta)$.

We will compare the numerators and the denominators of $FDR$ in (8) and $\delta_n$ in (1) respectively. Since the comparisons are almost identical, we discuss the comparison between the numerators only. We denote $E_\theta(\cdot)$ and $V_\theta(\cdot)$ as the conditional mean and variance given the true parameter $\theta$, and we denote $E(\cdot)$ and $V(\cdot)$ as the marginal mean and variance under the prior $G(\theta)$. Let $Y_i = I_{T_{ni} \in C, \theta_i \le \theta_0}$. Then given $\theta_1, \ldots, \theta_m$, $Y_i$ $(i = 1, \ldots, m)$ are independent Bernoulli random variables with mean values $\mu_i = \mu_i(\theta_i) = E_{\theta_i}(Y_i)$, and marginally $\mu_i$ are iid with expected value $A_n$ in (3). Let

$$D_m = \frac{1}{m} \sum_{i=1}^m I_{T_{ni} \in C, \theta_i \le \theta_0} - A_n = \frac{1}{m} \sum_{i=1}^m (Y_i - \mu_i) + \frac{1}{m} \sum_{i=1}^m (\mu_i - A_n).$$

Note that we assume that $\theta_1, \ldots, \theta_m$ are iid with a common CDF $G(\theta)$. The second term goes to 0 almost surely by the Strong Law of Large Numbers (SLLN) for identically distributed random variables. Note that for any given $\theta_1, \ldots, \theta_m$, $Y_1, \ldots, Y_m$ are independent but not iid, with $E_{\theta_i}(Y_i) = \mu_i$, $V_{\theta_i}(Y_i) = \mu_i(1 - \mu_i)$ and $\sum_{i=1}^{-\infty} i^{-2} V_{\theta_i}(Y_i) \leq \sum_{i=1}^{\infty} i^{-2} < \infty$. The first term also goes to 0 almost surely by a SLLN for independent but not iid random variables [15]. Therefore, $D_m$ goes to 0 almost surely. The comparison of denominators is handled similarly. Therefore, for almost all sequences $\theta_1, \theta_2, \ldots,$

$$\frac{\sum_{i=1}^m I_{T_{ni} \in C, \theta_i \leq \theta_0}}{(\sum_{i=1}^m I_{T_{ni} \in C}) \vee 1} \to \delta_n$$

as $m \to \infty$.

Since $\sum_{i=1}^m I_{T_{ni} \in C, \theta_i \leq \theta_0} \leq (\sum_{i=1}^m I_{T_{ni} \leq C}) \vee 1$, their ratio is uniformly integrable. And so, FDR as defined in (8) also converges to $\delta_n$ as $m \to \infty$ for almost all sequences $\theta_1, \theta_2, \ldots$.

This gives a pleasant, exact connection between our approach and the established indices formulated by the previous researchers. Of course, for fixed $m$, the frequentist FDR does not need to be close to our $\delta_n$.

## 3. Continuous one-parameter exponential family

Assume the density of the i.i.d. sample $X_1, \ldots, X_n$ is in the form of a one-parameter exponential family $f_\theta(x) = b(x)e^{\theta x - a(\theta)}$ for $x \in \mathcal{X} \subseteq R$, where the natural space $\Omega$ of $\theta$ is an interval of $R$ and $a(\theta) = \log \int_{\mathcal{X}} b(x)e^{\theta x} dx$. Without loss of generality, we can assume $\Omega$ is open so that one can write $\Omega = (\underline{\theta}, \bar{\theta})$ for $-\infty \leq \underline{\theta} < \bar{\theta} \leq \infty$. All derivatives of $a(\theta)$ exist at every $\theta \in \Omega$ and can be derived by formally differentiating under the integral sign ([4], p. 34). This implies that $a'(\theta) = E_\theta(X_1)$, $a''(\theta) = Var_\theta(X_1)$ for every $\theta \in \Omega$. Let us denote $\mu(\theta) = a'(\theta)$, $\sigma(\theta) = \sqrt{a''(\theta)}$, $\kappa_i(\theta) = a^{(i)}(\theta)$ and $\rho_i(\theta) = \kappa_i(\theta)/\sigma^i(\theta)$ for $i \geq 3$, where $a^{(i)}(\theta)$ represents the $i$-th derivative of $a(\theta)$. Then, $\mu(\theta)$, $\sigma(\theta)$, $\kappa_i(\theta)$ and $\rho_i(\theta)$ all exist and are continuous at every $\theta \in \Omega$ ([4], p. 36), and $\mu(\theta)$ is non-decreasing in $\theta$ since $a''(\theta) = \sigma^2(\theta) \geq 0$ for all $\theta$.

Let $\mu_0 = \mu(\theta_0)$, $\sigma_0 = \sigma(\theta_0)$, $\kappa_{i0} = \kappa_i(\theta_0)$ and $\rho_{i0} = \rho_i(\theta_0)$ for $i \geq 3$ and assume $\sigma_0 > 0$. The usual $\alpha$ $(0 < \alpha < 1)$ level UMP test ([13], p. 80) for the testing problem $H_0 : \theta \leq \theta_0$ vs $H_A : \theta \geq \theta_0$ rejects $H_0$ if $\bar{X} \in C$ where

$$(9) \qquad C = \{\bar{X} : \sqrt{n} \frac{\bar{X} - \mu_0}{\sigma_0} > k_{\theta_0, n}\},$$

and $k_{\theta_0, n}$ is determined from $P_{\theta_0}\{\sqrt{n}(\bar{X} - \mu_0)/\sigma_0 > k_{\theta_0, n}\} = \alpha$; $\lim_{n \to \infty} k_{\theta_0, n} = z_\alpha$. Let

$$(10) \qquad \tilde{\beta}_n(\theta) = P_\theta\left(\sqrt{n} \frac{\bar{X} - \mu_0}{\sigma_0} > k_{\theta_0, n}\right)$$

Then, using the transformation $x = \sigma_0 \sqrt{n}(\theta - \theta_0) - z_\alpha$ under the integral sign below, we have

$$(11) \qquad A_n = \int_{\underline{\theta}}^{\theta_0} \tilde{\beta}_n(\theta) g(\theta) d\theta = \frac{1}{\sigma_0 \sqrt{n}} \int_{\underline{x}}^{-z_\alpha} \tilde{\beta}_n(\theta_0 + \frac{x + z_\alpha}{\sigma_0 \sqrt{n}}) g(\theta_0 + \frac{x + z_\alpha}{\sigma_0 \sqrt{n}}) dx$$

and

$$(12) \qquad \tilde{A}_n = \frac{1}{\sigma_0\sqrt{n}} \int_{-z_\alpha}^{\bar{x}} [1 - \tilde{\beta}_n(\theta_0 + \frac{x + z_\alpha}{\sigma_0\sqrt{n}})]g(\theta_0 + \frac{x + z_\alpha}{\sigma_0\sqrt{n}})dx,$$

where $\underline{x} = \sigma_0\sqrt{n}(\underline{\theta} - \theta_0) - z_\alpha$ and $\bar{x} = \sigma_0\sqrt{n}(\bar{\theta} - \theta_0) - z_\alpha$.

Since for an interior parameter $\theta$ all moments of the exponential family exist and are continuous in $\theta$, we can find $\theta_1$ and $\theta_2$ satisfying $\bar{\theta} < \theta_1 < \theta_0$ and $\theta_0 < \theta_2 < \bar{\theta}$ such that for any $\theta \in [\theta_1, \theta_2]$, $\sigma^2(\theta)$, $\kappa_3(\theta)$, $\kappa_4(\theta)$, $\kappa_5(\theta)$, $g(\theta)$, $g'(\theta)$, $g''(\theta)$ and $g^{(3)}(\theta)$ are uniformly bounded in absolute values, and the minimum value of $\sigma^2(\theta)$ is a positive number. After we pick $\theta_1$ and $\theta_2$, we partition each of $A_n$ and $\tilde{A}_n$ into two parts so that one part is negligible in the expansion. Then, the rest of the work in the expansion is to find the coefficients of the second part.

To describe these partitions, we define $\theta_{1n} = \theta_0 + (\theta_1 - \theta_0)/n^{1/3}$, $\theta_{2n} = \theta_0 + (\theta_2 - \theta_0)/n^{1/3}$, $x_{1n} = \sigma_0\sqrt{n}(\theta_{1n} - \theta_0) - z_\alpha$ and $x_{2n} = \sigma_0\sqrt{n}(\theta_{2n} - \theta_0) - z_\alpha$. Let

$$(13) \qquad A_{n,\theta_{1n}} = \frac{1}{\sigma_0\sqrt{n}} \int_{x_{1n}}^{-z_\alpha} \tilde{\beta}_n(\theta_0 + \frac{x + z_\alpha}{\sigma_0\sqrt{n}})g(\theta_0 + \frac{x + z_\alpha}{\sigma_0\sqrt{n}})dx$$

$$(14) \qquad \underline{R}_{n,\theta_{1n}} = \frac{1}{\sigma_0\sqrt{n}} \int_{\underline{x}}^{x_{1n}} \tilde{\beta}_n(\theta_0 + \frac{x + z_\alpha}{\sigma_0\sqrt{n}})g(\theta_0 + \frac{x + z_\alpha}{\sigma_0\sqrt{n}})dx,$$

$$(15) \qquad \tilde{A}_{n,\theta_{2n}} = \frac{1}{\sigma_0\sqrt{n}} \int_{-z_\alpha}^{x_{2n}} [1 - \tilde{\beta}_n(\theta_0 + \frac{x + z_\alpha}{\sigma_0\sqrt{n}})]g(\theta_0 + \frac{x + z_\alpha}{\sigma_0\sqrt{n}})dx,$$

and

$$(16) \qquad \bar{R}_{n,\theta_{2n}} = \frac{1}{\sigma_0\sqrt{n}} \int_{x_{2n}}^{\bar{x}} [1 - \tilde{\beta}_n(\theta_0 + \frac{x + z_\alpha}{\sigma_0\sqrt{n}})]g(\theta_0 + \frac{x + z_\alpha}{\sigma_0\sqrt{n}})dx.$$

Then, $A_n = A_{n,\theta_{1n}} + \underline{R}_{n,\theta_{1n}}$ and $\tilde{A}_n = \tilde{A}_{n,\theta_{2n}} + \bar{R}_{n,\theta_{2n}}$. In the appendix, we show that for any $\ell > 0$, $\lim_{n\to\infty} n^\ell \underline{R}_{n,\theta_{1n}} = \lim_{n\to\infty} n^\ell \bar{R}_{n,\theta_{2n}} = 0$. Therefore, it is enough to compute the coefficients of the expansions for $A_{n,\theta_{1n}}$ and $\tilde{A}_{n,\theta_{2n}}$. Among the steps for expansions, the key step is to compute the expansions of $\tilde{\beta}_n(\theta_0 + (x + z_\alpha)/(\sigma_0\sqrt{n}))$ when $x \in [x_{1n}, -z_\alpha]$ and $1 - \tilde{\beta}_n(\theta_0 + (x + z_\alpha)/(\sigma_0\sqrt{n}))$ when $x \in [-z_\alpha, x_{2n}]$ under the integral sign, since the expansion of $g(\theta_0 + (x + z_\alpha)/(\sigma_0\sqrt{n}))$ in (13) and (15) is easily obtained as

$$(17) \qquad g(\theta_0 + \frac{x + z_\alpha}{\sigma_0\sqrt{n}}) = g(\theta_0) + g'(\theta_0)\frac{x + z_\alpha}{\sigma_0\sqrt{n}} + \frac{g''(\theta_0)}{2}\frac{(x + z_\alpha)^2}{\sigma_0^2 n} + O(n^{-2}).$$

After a lengthy calculation, we have

$$(18) \qquad \begin{aligned} A_{n,\theta_{1n}} = \frac{1}{\sigma_0\sqrt{n}} \int_{x_{1n}}^{-z_\alpha} &[\Phi(x) + \frac{\phi(x)g_1(x)}{\sqrt{n}} + \frac{\phi(x)g_2(x)}{n}] \\ &\times [g(\theta_0) + g'(\theta_0)\frac{x + z_\alpha}{\sigma_0\sqrt{n}} + \frac{g''(\theta_0)}{2}\frac{(x + z_\alpha)^2}{\sigma_0^2 n}]dx + O(n^{-2}), \end{aligned}$$

and

$$(19) \qquad \begin{aligned} \tilde{A}_{n,\theta_{2n}} = \frac{1}{\sigma_0\sqrt{n}} \int_{-z_\alpha}^{x_{2n}} &[1 - \Phi(x) - \frac{\phi(x)g_1(x)}{\sqrt{n}} - \frac{\phi(x)g_2(x)}{n}] \\ &\times [g(\theta_0) + g'(\theta_0)\frac{x + z_\alpha}{\sigma_0\sqrt{n}} + \frac{g''(\theta_0)}{2}\frac{(x + z_\alpha)^2}{\sigma_0^2 n}]dx + O(n^{-2}). \end{aligned}$$

where

$$(20) \qquad g_1(x) = \frac{\rho_{30}}{6} x^2 + \frac{z_\alpha \rho_{30}}{2} x + \frac{z_\alpha^2 \rho_{30}}{3}$$

and

$$
\begin{aligned}
(21) \quad g_2(x) = {} & \frac{\rho_{30}^2}{72} x^5 - \frac{z_\alpha \rho_{30}^2}{12} x^4 + \left(\frac{\rho_{40}}{8} - \frac{13 z_\alpha^2 \rho_{30}^2}{72} - \frac{7 \rho_{30}^2}{24}\right) x^3 + \left(\frac{z_\alpha \rho_{40}}{6} - \frac{z_\alpha^3 \rho_{30}^2}{6}\right. \\
& \left. - \frac{z_\alpha \rho_{30}^2}{12}\right) x^2 + \left[\left(\frac{z_\alpha^2}{4} - \frac{7}{24}\right) \rho_{40} - \frac{z_\alpha^4 \rho_{30}^2}{18} - \frac{13 z_\alpha^2 \rho_{30}^2}{72} + \frac{4 \rho_{30}^2}{9}\right] x \\
& + \left[\left(\frac{z_\alpha^3}{8} - \frac{z_\alpha}{24}\right) \rho_{40} - \left(\frac{z_\alpha^3}{9} - \frac{z_\alpha}{36}\right) \rho_{30}^2\right].
\end{aligned}
$$

The expressions for $g_1(x)$ and $g_2(x)$ are derived in the Appendix; the derivation of these two formulae forms the dominant part of the penultimate expression and involves the use of Cornish–Fisher as well as Edgeworth expansions.

On using (18), (19), (20) and (21), we have the following expansions

$$(22) \qquad A_{n,\theta_{1n}} = \frac{a_1}{\sqrt{n}} + \frac{a_2}{n} + \frac{a_3}{n^{3/2}} + O(n^{-2}),$$

where

$$
\begin{aligned}
(23) \quad a_1 = {} & \frac{g(\theta_0)}{\sigma_0} [\phi(z_\alpha) - \alpha z_\alpha], \\
a_2 = {} & \frac{\rho_{30} g(\theta_0)}{6\sigma_0} [\alpha + 2\alpha z_\alpha^2 - 2 z_\alpha \phi(z_\alpha)] - \frac{g'(\theta_0)}{2\sigma_0^2} [\alpha(z_\alpha^2 + 1) - z_\alpha \phi(z_\alpha)] \\
a_3 = {} & \frac{g''(\theta_0)}{6\sigma_0^3} [(z_\alpha^2 + 2)\phi(z_\alpha) - \alpha(z_\alpha^3 + 3 z_\alpha)] + \frac{g'(\theta_0)}{\sigma_0^2} \Big[\frac{\alpha \rho_{30}}{3}(z_\alpha^3 + 2 z_\alpha) \\
& - \frac{\rho_{30}}{3}(z_\alpha^2 + 1)\phi(z_\alpha)\Big] + \frac{g(\theta_0)}{\sigma_0} \Big[\left(-\frac{z_\alpha^4 \rho_{30}^2}{36} + \frac{4 z_\alpha^2 \rho_{30}^2}{9} + \frac{\rho_{30}^2}{36}\right. \\
& \left. - \frac{5 z_\alpha^2 \rho_{40}}{24} + \frac{\rho_{40}}{24}\right)\phi(z_\alpha) + \alpha\left(-\frac{5 z_\alpha^3 \rho_{30}^2}{18} - \frac{11 z_\alpha \rho_{30}^2}{36} + \frac{z_\alpha^3 \rho_{40}}{8} + \frac{z_\alpha \rho_{40}}{8}\right)\Big].
\end{aligned}
$$

Similarly,

$$(24) \qquad \tilde{A}_{n,\theta_{2n}} = \frac{\tilde{a}_1}{\sqrt{n}} + \frac{\tilde{a}_2}{n} + \frac{\tilde{a}_3}{n^{3/2}} + O(n^{-2}),$$

where $\tilde{a}_1 = [g(\theta_0)/\sigma_0][\phi(z_\alpha) + (1-\alpha) z_\alpha]$, $\tilde{a}_2 = [g'(\theta_0)/(2\sigma_0^2)][(1-\alpha)(z_\alpha^2 + 1) + z_\alpha \phi(z_\alpha)] - [\rho_{30} g(\theta_0)/(6\sigma_0)][(1-\alpha)(1 + 2 z_\alpha^2) + 2 z_\alpha \phi(z_\alpha)]$, $\tilde{a}_3 = [g''(\theta_0)/(6\sigma_0^3)][(z_\alpha^2 + 2)\phi(z_\alpha) + (1-\alpha)(z_\alpha^3 + 3 z_\alpha)] - [g'(\theta_0)\rho_{30}/(3\sigma_0^2)][(z_\alpha^2 + 1)\phi(z_\alpha) + (1-\alpha)(z_\alpha^3 + 2 z_\alpha)] + [g(\theta_0)/\sigma_0][\phi(z_\alpha)(-z_\alpha^4 \rho_{30}^2/36 + 4 z_\alpha^2 \rho_{30}^2/9 + \rho_{30}^2/36 - 5 z_\alpha^2 \rho_{40}/24 + \rho_{40}/24) - (1-\alpha)(-5 z_\alpha^3 \rho_{30}^2/18 - 11 z_\alpha \rho_{30}^2/36 + z_\alpha^3 \rho_{40}/8 + z_\alpha \rho_{40}/8)]$. The details of the expansions for $A_{n,\theta_{1n}}$ and $\tilde{A}_{n,\theta_{2n}}$ are given in the Appendix. Because the remainders $\underline{R}_{n,\theta_{1n}}$ and $\bar{R}_{n,\theta_{2n}}$ are of smaller order than $n^{-2}$ as we commented before, the expansions in (22) and (24) are the expansions for $A_n$ and $\tilde{A}_n$ in (3) and (4) respectively.

The expansions of $\delta_n$ and $\epsilon_n$ in (1) and (2) can now be obtained by letting $\lambda = \int_{\underline{\theta}}^{\theta_0} g(\theta) d\theta$, $b_1 = \tilde{a}_1 - a_1$, $b_2 = \tilde{a}_2 - a_2$ and $b_3 = \tilde{a}_3 - a_3$ in (7).

### *3.1. Examples*

**Example 1.** Let $X_1, \ldots, X_n$ be i.i.d. $N(\theta, 1)$. Since $\theta$ is a location parameter, there is no loss of generality in letting $\theta_0 = 0$. Thus consider testing $H_0 : \theta \leq 0$ vs $H_1 : \theta > 0$. Clearly, we have $\mu(\theta) = \theta$, $\sigma(\theta) = 1$ and $\rho_i(\theta) = \kappa_i(\theta) = 0$ for all $i \geq 3$.

The $\alpha$ $(0 < \alpha < 1)$ level UMP test rejects $H_0$ if $\sqrt{n}\bar{X} > z_\alpha$. For a continuously three times differentiable prior $g(\theta)$ for $\theta$, one can simply plug the values of $\mu_0 = 0$, $\sigma_0 = 1$, $\rho_{30} = \rho_{40} = 0$ into (23) and the coefficients of the expansion in (24) to get the coefficients $a_1 = g(0)[\Phi(z_\alpha) - \alpha z_\alpha]$, $a_2 = -g'(0)[\alpha(z_\alpha^2 - 1) - z_\alpha\phi(z_\alpha)]$, $a_3 = g''(0)[(z_\alpha + 2)\phi(z_\alpha) - \alpha(z_\alpha^3 + 3z_\alpha)]/6$, $\tilde{a}_1 = g(0)[\phi(z_\alpha) + \phi(z_\alpha)]$, $\tilde{a}_2 = g'(0)[(1-\alpha)(z_\alpha^2 + 1) + z_\alpha\phi(z_\alpha)]/2$, $\tilde{a}_3 = g''(0)[(z_\alpha + 2)\phi(z_\alpha) + (1-\alpha)(z_\alpha^3 + 3z_\alpha)]/6$. Substituting $a_1$, $a_2$, $a_3$, $\tilde{a}_1$, $\tilde{a}_2$ and $\tilde{a}_3$ into (7), one derives the expansions of $\delta_n$ and $\epsilon_n$ as given by (1) and (2) respectively.

If the prior density function is also assumed to be symmetric, then $\lambda = 1/2$ and $g'(0) = 0$. In this case, the coefficients of the expansion of $\delta_n$ in (1) are given explicitly as follows: $c_1 = 2g(0)[\phi(z_\alpha) - \alpha z_\alpha]$, $c_2 = 4z_\alpha[g(0)]^2[\phi(z_\alpha) - \alpha z_\alpha]$, $c_3 = 2\phi(z_\alpha)\{4z_\alpha^2[g(0)]^3 + g''(0)(z_\alpha^2 + 2)/6\} - \alpha\{g''(0)(z_\alpha^3 + 3z_\alpha)/3 + 8z_\alpha^3[g(0)]^3\}$, and the coefficients of the expansions of $\epsilon_n$ in (2) are as $d_1 = 2g(0)[(1-\alpha)z_\alpha + \phi(z_\alpha)]$, $d_2 = -4z_\alpha[g(0)]^2[(1-\alpha)z_\alpha + \phi(z_\alpha)]$, $d_3 = 2\phi(z_\alpha)\{4z_\alpha^2[g(0)]^3 + g''(0)(z_\alpha^2 + 2)/6\} + (1-\alpha)\{g''(0)(z_\alpha^3 + 3z_\alpha)/3 + 8z_\alpha^3[g(0)]^3\}$.

Two specific prior distributions for $\theta$ are now considered for numerical illustration. In the first one we choose $\theta \sim N(0, \tau^2)$ and in the second example we choose $\theta/\tau \sim t_m$, where $\tau$ is a scale parameter. Clearly $g^{(3)}(\theta)$ is continuous in $\theta$ in both cases.

If $g(\theta)$ is the density of $\theta$ when $\theta \sim N(0, \tau^2)$, then $\lambda = 1/2$, $g(0) = 1/[\sqrt{2\pi}\tau]$, $g'(0) = 0$ and $g''(0) = -1/[\sqrt{2\pi}\tau^3]$.

We calculated the numerical values of $c_1$, $c_2$, $c_3$, $d_1$, $d_2$ and $d_3$ as functions of $\alpha$ when $\theta \sim N(0, 1)$. We note that $c_1$ is a monotone increasing function and $d_1$ is also a monotone decreasing function of $\alpha$. However, $c_2$, $d_2$ and $c_3$, $d_3$ are not monotone and in fact, $d_2$ is decreasing when $\alpha$ is close to 1 (not shown), $c_3$ also takes negative values and $d_3$ takes positive values for larger values of $\alpha$.

If $g(\theta)$ is the density of $\theta$ when $\theta/\tau \sim t_m$, then $\lambda = 1/2$, $g'(0) = 0$, $g(0) = \Gamma(\frac{m+1}{2})/[\tau\sqrt{m\pi}\Gamma(\frac{m}{2})]$ and $g''(0) = -\Gamma(\frac{m+3}{2})/[\tau\sqrt{m\pi}\Gamma(\frac{m+2}{2})]$. Putting those values into the corresponding expressions, we get the coefficients $c_1, c_2, c_3$ and $d_1, d_2, d_3$ of the expansions of $\delta_n$ and $\epsilon_n$. When $m = 1$, the results are exactly the same as the Cauchy prior for $\theta$.

Numerical results very similar to the normal prior are seen for the Cauchy case. From Figure 1, we see that for each of the normal and the Cauchy prior, *only about 1% of those null hypotheses a frequentist rejects with a p-value of less than 5% are true*. Indeed quite typically, $\delta_n < \alpha$ for even very small values of $n$. This is discussed in greater detail in Section 4.5. This finding seems to be quite interesting.

The true values of $\delta_n$ and $\epsilon_n$ are computed by taking an average of the lower and the upper Riemann sums in $A_n$, $\tilde{A}_n$, $B_n$ and $\tilde{B}_n$ with the exact formulae for the standard normal pdf. The accuracy of the expansion for $\delta_n$ is remarkable, as can be seen in Figure 1. Even for $n = 4$, the true value of $\delta_n$ is almost identical to the expansion in (1). The accuracy of the expansion for $\epsilon_n$ is very good (even if it is not as good as that for $\delta_n$). For $n = 20$, the true value of $\epsilon_n$ is almost identical to the expansion in (2).

**Example 2.** Let $X_1, \cdots, X_n$ be iid $Exp(\theta)$, with density $f_\theta(x) = \theta e^{-\theta x}$ if $x > 0$. Clearly, $\mu(\theta) = 1/\theta$, $\sigma^2(\theta) = 1/\theta^2$, $\rho_3(\theta) = 2$ and $\rho_4(\theta) = 6$. Let $\tilde{\theta} = -\theta$. Then, one can write the density of $X_1$ in the standard form of the exponential family as $f_{\tilde{\theta}}(x) = e^{\tilde{\theta} x + \log|\tilde{\theta}|}$. The natural parameter space of $\tilde{\theta}$ is $\Omega = (-\infty, 0)$. If $g(\theta)$ is a prior density for $\theta$ on $(0, \infty)$, then $g(-\tilde{\theta})$ is a prior density for $\tilde{\theta}$ on $(-\infty, 0)$. Since $\theta$ is a scale parameter, it is enough to look at the case $\tilde{\theta}_0 = -1$. In terms of $\theta$, therefore the problem considered is to test $H_0 : \theta \geq 1$ vs $H_1 : \theta < 1$. The $\alpha$ $(0 < \alpha < 1)$ level UMP test for this problem rejects $H_0$ if $\bar{X} > \Gamma_{\alpha,n,n}$, where $\Gamma_{\alpha,r,s}$ is the upper $\alpha$ quantile of the Gamma distribution with parameters $r$ and $s$. If $g(\theta)$ is continuous and three time differentiable, then we can simply put the values $\mu_0 = 1$, $\sigma_0 = 1$, $\rho_{30} = 2$, $\rho_{40} = 6$, and $\lambda = \int_0^1 g(\theta)d\theta$ into (23) and the coefficients of the expansion in (24) to get the coefficients $a_1$, $a_2$, $a_3$, $\tilde{a}_1$, $\tilde{a}_2$ and $\tilde{a}_3$, and then get the expansions of $\delta_n$ and $\epsilon_n$ in (1) and (2) respectively.

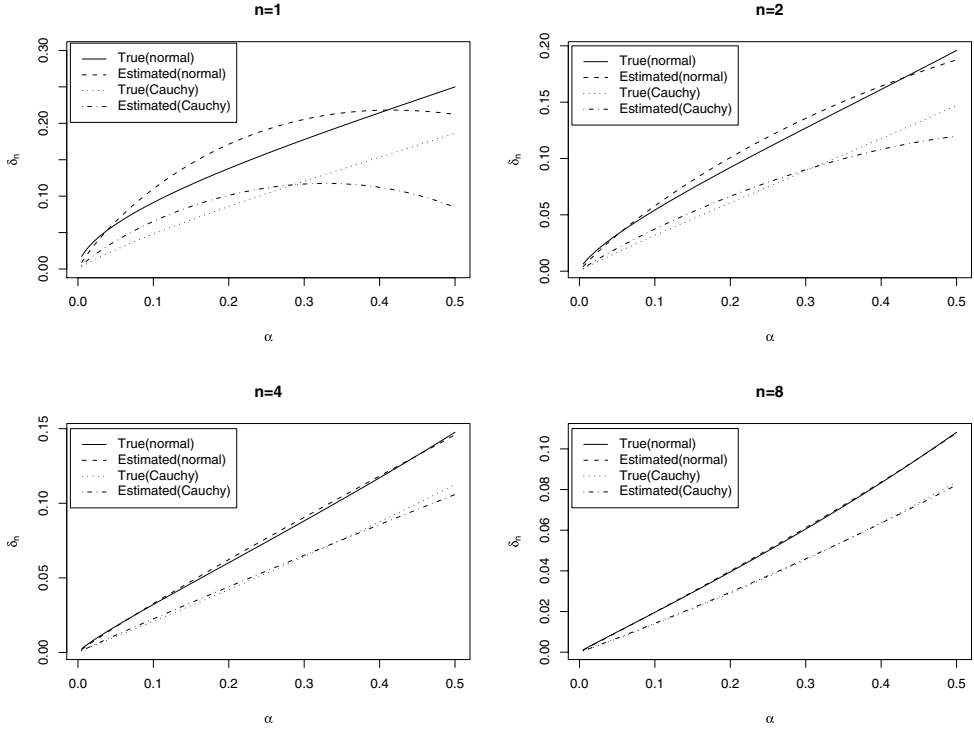

FIG 1. *True and estimated values of $\delta_n$ as functions of $\alpha$ for the standard normal prior and the Cauchy prior.*

Two priors are to be considered in this example. The first one is the Gamma prior with prior density $g(\theta) = s^r \theta^{r-1} e^{-s\theta} / \Gamma(r)$, where $r$ and $s$ are known constants. It would be natural to have the mode of $g(\theta)$ at 1, that is $s = r - 1$. In this case, $g'(1) = 0$, $g(1) = (r-1)^r e^{-(r-1)} / \Gamma(r)$ and $g''(1) = -(r-1)^{r+1} e^{-(r-1)} / \Gamma(r)$.

Next, consider the $F$ prior with degrees of freedom $2r$ and $2s$ for $\theta/\tau$ for a fixed $\tau > 0$. Then, the prior density for $\theta$ is $g(\theta) = \frac{\Gamma(r+s)}{\Gamma(r)\Gamma(s)} \frac{r}{s\tau} (\frac{r\theta}{s\tau})^{r-1} (1 + \frac{r\theta}{s\tau})^{-(r+s)}$. To make the mode of $g(\theta)$ equal to 1, we have to choose $\tau = r(s+1)/[s(r-1)]$. Then $g'(1) = 0$, $g(1) = \frac{\Gamma(r+s)}{\Gamma(r)\Gamma(s)} (\frac{r-1}{s+1})^r (1 + \frac{r-1}{s+1})^{-(r+s)}$, and $g''(1) = -\frac{\Gamma(r+s)}{\Gamma(r)\Gamma(s)} (\frac{r-1}{s+1})^{r+1} (r+s)(1 + \frac{r-1}{s+1})^{-(r+s+2)}$.

Exact and estimated values of $\delta_n$ are plotted in Figure 3. At $n = 20$, the expansion is clearly extremely accurate and as in example 1, we see that the false

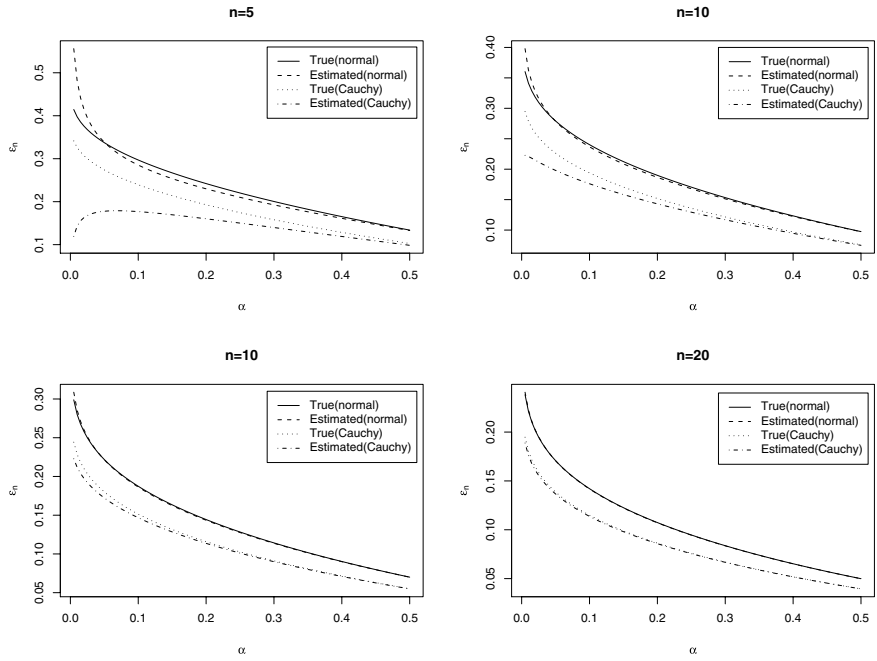

FIG 2. *True and estimated values of $\epsilon_n$ as functions of $\alpha$ for the standard normal prior and the Cauchy prior.*

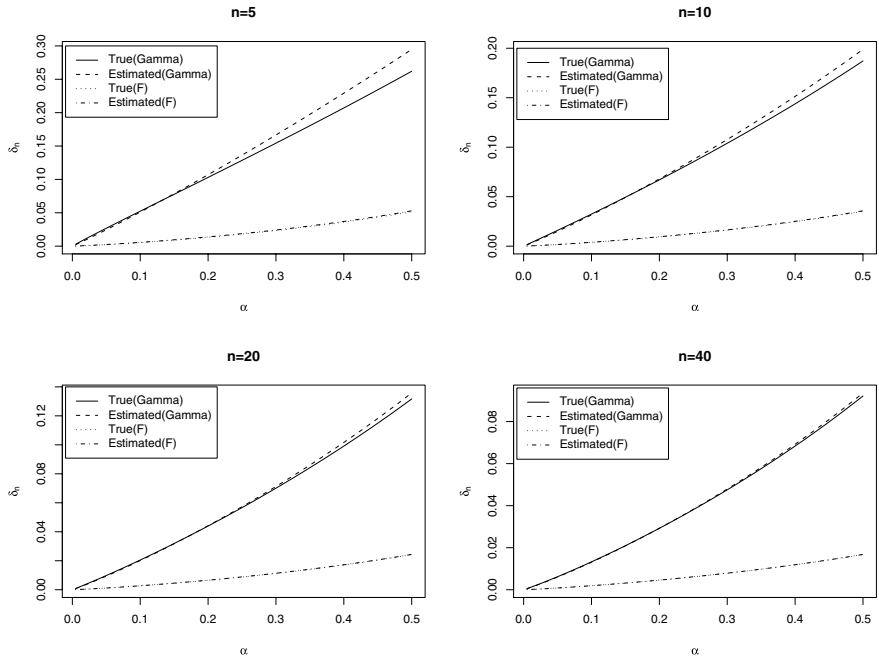

FIG 3. *True and estimated values of $\delta_n$ as functions of $\alpha$ under $\Gamma(2, 1)$ and $F(4, 4)$ priors for $\theta$ when $X \sim Exp(\theta)$.*

discovery rate $\delta_n$ is very small even for $n = 10$.

### 3.2. The frequentist is more prone to type II error

Consider the two Bayesian error rates

$$\delta_n = P(H_0 \mid \text{Frequentist rejects } H_0)$$

and

$$\epsilon_n = P(H_1 \mid \text{Frequentist accepts } H_0).$$

Is there an inequality between $\delta_n$ and $\epsilon_n$? Rather interestingly, when $\theta$ is the normal mean and the testing problem is $H_0 : \theta \leq 0$ versus $H_1 : \theta > 0$, there is an approximate inequality in the sense that if we consider the respective coefficients $c_1$ and $d_1$ of the $1/\sqrt{n}$ term, then for any symmetric prior (because then $g'(0) = 0$ and $\lambda = 1 - \lambda = 1/2$), we have

$$c_1 = 2g(0)[\phi(z_\alpha) - \alpha z_\alpha] < d_1 = 2g(0)[(1-\alpha)z_\alpha + \phi(z_\alpha)]$$

for any $\alpha < 1/2$. It is interesting that this inequality holds regardless of the exact choice of $g(\cdot)$ and the value of $\alpha$, as long as $\alpha < 1/2$. Thus, to the first order, the frequentist is less prone to type I error. Even the exact values of $\delta_n$ and $\epsilon_n$ satisfy this inequality, unless $\alpha$ is small, as can be seen, for example from a scrutiny of Figures 1 and 2. This would suggest that a frequentist needs to be more mindful of premature acceptance of $H_0$ rather than its premature rejection in the composite one sided problem. This is in contrast to the conclusion reached in Berger and Sellke [2] under their formulation.

## 4. General location parameter case

As we mentioned in Section 1, the quantities $\delta_n, \epsilon_n$ depend on the choice of the test statistic. For location parameter problems, in general there is no reason to use the sample mean as the test statistic. For many non-normal location parameter densities, such as the double exponential, it is more natural to use the sample median as the test statistic.

Assume the density of the i.i.d. sample $X_1, \ldots, X_n$ is $f(x - \theta)$ where the median of $f(\cdot)$ is 0, and assume $f(0) > 0$. Then an asymptotic size $\alpha$ test for

$$H_0 : \theta \leq 0 \text{ vs } H_1 : \theta > 0$$

rejects $H_0$ if $\sqrt{n}T_n > z_\alpha/[2f(0)]$, where $T_n = X_{([\frac{n}{2}]+1)}$ is the sample median ([8], p. 89), since $\sqrt{n}(T_n - \theta) \overset{L}{\Rightarrow} N(0, 1/[4f^2(0)])$. We will derive the coefficients $c_1, c_2, c_3$ in (1) and $d_1, d_2, d_3$ in (2) given the prior density $g(\theta)$ for $\theta$. We assume again that $g(\theta)$ is three times differentiable with a bounded absolute third derivative.

### 4.1. Expansion of type I error and type II error

To obtain the coefficients of the expansions of $\delta_n$ in (1) and $\epsilon_n$ in (2), we have to expand the $A_n$ and $\tilde{A}_n$ given by (3) and (4). Of these,

$$(25) \quad A_n = P(\theta \leq 0, \sqrt{n}T_n > \frac{z_\alpha}{2f(0)}) = \frac{1}{\sqrt{n}} \int_{-\infty}^{0} \{1 - F_n[z_\alpha - 2xf(0)]\} g(\frac{x}{\sqrt{n}}) dx$$

where $F_n$ is the CDF of $2f(0)\sqrt{n}(T_n - \theta)$ if the true median is $\theta$. Reiss [16] gives the expansion of $F_n$ as

(26) $$F_n(t) = \Phi(t) + \frac{\phi(t)}{\sqrt{n}}R_1(t) + \frac{\phi(t)}{n}R_2(t) + r_{t,n},$$

where, with $\{x\}$ denoting the fractional part of a real $x$, $R_1(t) = f_{11}t^2 + f_{12}$, $f_{11} = f'(0)/[4f^2(0)]$, $f_{12} = -(1 - 2\{\frac{n}{2}\})$, and $R_2(t) = f_{21}t^5 + f_{22}t^3 + f_{23}t$, where $f_{21} = -[f'(0)/f^2(0)]^2/32$, $f_{22} = 1/4 + (1/2 - \{\frac{n}{2}\})[f'(0)/(2f^2(0))] + f''(0)/[24f^3(0)]$, $f_{23} = 1/4 - (1 - 2\{\frac{n}{2}\})^2/2$. The error term $r_{1,t,n}$ can be written as $r_{t,n} = \phi(t)R_3(t)/n^{3/2} + O(n^{-2})$, where $R_3(t)$ is a polynomial.

By letting $y = 2xf(0) - z_\alpha$ in (25), we have

(27) $$A_n = \frac{1}{2f(0)\sqrt{n}} \int_{-\infty}^{-z_\alpha} \{\Phi(y) - \frac{\phi(y)}{\sqrt{n}}R_1(-y) - \frac{\phi(y)}{n}R_2(-y) - r_{-y,n}\}$$
$$\times [g(0) + g'(0)\frac{y + z_\alpha}{2f(0)\sqrt{n}} + \frac{g''(0)}{2}\frac{(y + z_\alpha)^2}{4f^2(0)n} + \frac{(y + z_\alpha)^3}{48f^3(0)n^{3/2}}g^{(3)}(y^*)]dy,$$

where $y^*$ is between 0 and $(y + z_\alpha)/[2f(0)\sqrt{n}]$.

Hence, assuming $\sup_\theta |g^{(3)}(\theta)| < \infty$, on exact integration of each product of functions in (27) and on collapsing the terms, we get

(28) $$A_n = \frac{a_1}{\sqrt{n}} + \frac{a_2}{n} + \frac{a_3}{n^{3/2}} + O(n^{-2}),$$

where

(29) $$a_1 = \frac{g(0)}{2f(0)}[\phi(z_\alpha) - \alpha z_\alpha],$$

(30) $$a_2 = \frac{g'(0)}{8f^2(0)}[z_\alpha\phi(z_\alpha) - \alpha(z_\alpha^2 + 1)] - \frac{g(0)}{2f(0)}\{f_{11}[z_\alpha\phi(z_\alpha) + \alpha] + f_{12}\alpha\}$$

and

(31)
$$a_3 = \frac{g''(0)}{48f^3(0)}[(z_\alpha^2 + 2)\phi(z_\alpha) - \alpha(z_\alpha^3 + 3z_\alpha)]$$
$$- \frac{g'(0)}{4f^2(0)}\{f_{11}[\alpha z_\alpha - 2\phi(z_\alpha)] + f_{12}[\alpha z_\alpha - \phi(z_\alpha)]\}$$
$$- \frac{g(0)}{2f(0)}\{f_{21}[(z_\alpha^4 + 4z_\alpha^2 + 8)\phi(z_\alpha)] + f_{22}[(z_\alpha^2 + 2)\phi(z_\alpha)] + f_{23}\phi(z_\alpha)\}.$$

We claim the error term in (28) is $O(n^{-2})$. To prove this, we need to look at its exact form, namely,

$$-\frac{O(n^{-2})}{2f(0)} \int_{-\infty}^{-z_\alpha} \phi(y)R_3(-y)g(\theta_0 + \frac{y + z_\alpha}{2f(0)\sqrt{n}})dy + O(n^{-2})\int_{-\infty}^{0} g(\theta_0 + y)dy.$$

Since $g(\theta)$ is absolutely uniformly bounded, the first term above is bounded by $O(n^{-2})$. The second term is $O(n^{-2})$ obviously. This shows that the error term in (28) is $O(n^{-2})$.

As regards $\tilde{A}_n$ given by (4), one can similarly obtain

$$
\begin{aligned}
\tilde{A}_n = P(\theta > 0, \sqrt{T}_n \leq \frac{z_\alpha}{2f(0)}) &= \frac{1}{\sqrt{n}} \int_0^\infty F_n[z_\alpha - 2f(0)x]g(\frac{x}{\sqrt{n}})dx \\
(32) \qquad\qquad &= \frac{\tilde{a}_1}{\sqrt{n}} + \frac{\tilde{a}_2}{n} + \frac{\tilde{a}_3}{n^{\frac{3}{2}}} + O(n^{-2}),
\end{aligned}
$$

where $y^*$ is between 0 and $(z_\alpha - y)/[2f(0)\sqrt{n}]$, $\tilde{a}_1 = [g(0)/(2f(0))][(1-\alpha)z_\alpha + \phi(z_\alpha)]$, $\tilde{a}_2 = [g'(0)/(8f^2(0))][(1-\alpha)(z_\alpha^2 + 1) + z_\alpha \phi(z_\alpha)] + [g(0)/(2f(0))]\{f_{11}[(1-\alpha) - z_\alpha \phi(z_\alpha)] + f_{12}(1-\alpha)\}$, $\tilde{a}_3 = [g''(0)/(48f^3(0))][(z_\alpha^2 + 2)\phi(z_\alpha) + (1-\alpha)(z_\alpha^3 + 3z_\alpha)] + [g'(0)/(4f^2(0))]\{f_{11}[(1-\alpha)z_\alpha + 2\phi(z_\alpha)] + f_{12}[(1-\alpha)z_\alpha + \phi(z_\alpha)]\} - [g(0)/(2f(0))] \times \{f_{21}[(z_\alpha^4 + 4z_\alpha^2 + 8)\phi(z_\alpha)] + f_{22}[(z_\alpha^2 + 2)\phi(z_\alpha)] + f_{23}\phi(z_\alpha)\}$. The error term in (32) is still $O(n^{-2})$ and this proof is omitted.

Therefore, we have the the expansions of $B_n$ given by (5) $B_n = \lambda - b_1/\sqrt{n} - b_2/n - b_3/n^{3/2} + O(n^{-2})$ where $\lambda = \int_0^\infty g(\theta)d\theta$ as before, $b_1 = \tilde{a}_1 - a_1 = z_\alpha g(0)/[2f(0)]$, $b_2 = \tilde{a}_2 - a_2 = g'(0)(z_\alpha^2 + 1)/[8f^2(0)] + g(0)(f_{11} + f_{12})/[2f(0)]$, $b_3 = \tilde{a}_3 - a_3 = g''(0)(z_\alpha^3 + 3z_\alpha)/[48f^3(0)] + z_\alpha g'(0)(f_{11} + f_{12})/[4f^2(0)]$. Substituting $a_1$, $a_2$, $a_3$, $\tilde{a}_1$, $\tilde{a}_2$, $\tilde{a}_3$, $b_1$, $b_2$ and $b_3$ into (7), we get the expansions of $\delta_n$ and $\epsilon_n$ for the general location parameter case given by (1) and (2).

## 4.2. Testing with mean vs. testing with median

Suppose $X_1, \ldots, X_n$ are i.i.d. observations from a $N(\theta, 1)$ density and the statistician tests $H_0 : \theta \leq 0$ vs. $H_1 : \theta > 0$ by using either the sample mean $\bar{X}$ or the median $T_n$. It is natural to ask the choice of which statistic makes him more vulnerable to false discoveries. We can look at both false discovery rates $\delta_n$ and $\epsilon_n$ to make this comparison, but we will do so only for the type I error rate $\delta_n$ here.

We assume for algebraic simplicity that $g$ is symmetric, and so $g'(0) = 0$ and $\lambda = 1/2$. Also, to keep track of the two statistics, we will denote the coefficients $c_1, c_2$ by $c_{1,\bar{X}}$, $c_{1,T_n}$, $c_{2,\bar{X}}$ and $c_{2,T_n}$ respectively. Then from our expansions in section 3.1 and section 4.1, it follows that

$$
c_{1,T_n} - c_{1,\bar{X}} = g(0)(\phi(z_\alpha) - \alpha z_\alpha)(\sqrt{2\pi} - 2) = a(\text{say}),
$$

and

$$
\begin{aligned}
c_{2,T_n} - c_{2,\bar{X}} &= g^2(0)z_\alpha(\phi(z_\alpha) - \alpha z_\alpha)(2\pi - 4) - g(0)\sqrt{2\pi}f_{12}\alpha \\
&\geq g^2(0)z_\alpha(\phi(z_\alpha) - \alpha z_\alpha)(2\pi - 4) = b(\text{say}) \text{ as } f_{12} \leq 0.
\end{aligned}
$$

Hence, there exist positive constants $a, b$ such that $\liminf_{n\to\infty} \sqrt{n}(\sqrt{n}(\delta_{n,T_n} - \delta_{n,\bar{X}}) - a) \geq b$, i.e., the statistician is *more* vulnerable to a type I false discovery by using the sample median as his test statistic. Now, of course, as a point estimator, $T_n$ is less efficient than the mean $\bar{X}$ in the normal case. Thus, the statistician is more vulnerable to a false discovery if he uses the less efficient point estimator as his test statistic. We find this neat connection between efficiency in estimation and false discovery rates in testing to be interesting. Of course, similar connections are well known in the literature on Pitman efficiencies of tests; see, e.g., van der Vaart ([24], p. 201).

### *4.3. Examples*

In this subsection, we are going to study the exact values and the expansions for $\delta_n$ and $\epsilon_n$ in two examples. One example is $f(x) = \phi(x)$ and $g(\theta) = \phi(\theta)$; for the other example, $f$ and $g$ are both densities of the standard Cauchy. We will refer to them as normal-normal and Cauchy-Cauchy for convenience of reference. The purpose of the first example is comparison with the normal-normal case when the test statistic was the sample mean (Example 2 in Section 3); the second example is an independent natural example.

For exact numerical evaluation of $\delta_n$ and $\epsilon_n$, the following formulae are necessary. The pdf of the standardized median $2f(0)\sqrt{n}(T_n - \theta)$ is

$$(33) \quad f_n(t) = \frac{\sqrt{n}\binom{n-1}{[\frac{n}{2}]-1}}{2f(0)} f(\frac{t}{2f(0)\sqrt{n}}) F^{[\frac{n}{2}]-1}(\frac{t}{2f(0)\sqrt{n}})(1 - F(\frac{t}{2f(0)\sqrt{n}}))^{n-[\frac{n}{2}]}.$$

We are now ready to present our examples.

**Example 3.** Suppose $X_1, X_2, \ldots, X_n$ are i.i.d. $N(\theta, 1)$ and $g(\theta) = \phi(\underline{\theta})$. Then, $g(0) = f(0) = 1/\sqrt{2\pi}$, $g'(0) = f'(0) = 0$ and $g''(0) = f''(0) = -1/\sqrt{2\pi}$. Then, we have $f_{11} = 0$, $f_{12} = -(1 - 2\{\frac{n}{2}\})$, $f_{21} = 0$, $f_{22} = 1/4 - \pi/12$ and $f_{23} = 1/4 - (1/2 - \{\frac{n}{2}\})^2$. Plugging these values for $f_{11}, f_{12}, f_{21}, f_{22}, f_{23}$ into (29), (30), (31) and (7), we obtain the expansions for $\delta_n$, and similarly for $\epsilon_n$ in the normal-normal case.

Next we consider the Cauchy-Cauchy case, i.e., $X_1, \ldots, X_n$ are i.i.d. with density function $f(x) = 1/\{\pi[1 + (x - \theta)^2]\}$ and $g(\theta) = 1/[\pi(1 + \theta^2)]$. Then, $f(0) = 1/\pi$, $f'(0) = 0$ and $f''(0) = -2/\pi$. Therefore, $f_{11} = 0$, $f_{12} = -(1 - 2\{\frac{n}{2}\})$, $f_{21} = 0$, $f_{22} = 1/4 - \pi^2/12$, and $f_{23} = 1/4 - (1/2 - \{\frac{n}{2}\})^2$. Plugging these values for $f_{11}, f_{12}, f_{21}, f_{22}, f_{23}$ in (29), (30), (31), we obtain the expansions for $\delta_n$, and similarly for $\epsilon_n$ in the Cauchy-Cauchy case.

The true and estimated values of $\delta_n$ for selected $n$ are given in Figure 4 and Figure 5. As before, the true values of $\delta_n$ and $\epsilon$ are computed by taking an average of the lower and the upper Riemann sums in $A_n$, $\tilde{A}_n$, $B_n$ and $\tilde{B}_n$ with the exact formulae for $f_n$ as in (33). It can be seen that the two values are almost identical when $n = 30$. By comparison with Figure 1, we see that the expansion for the median is not as precise as the expansion for the sample mean.

The most important thing we learn is how small $\delta_n$ is for very moderate values of $n$. For example, in Figure 4, $\delta_n$ is only about 0.01 if $\alpha = 0.05$, when $n = 20$. Again we see that even though we have changed the test statistic to the median, the frequentist's false discovery rate is very small and, in particular, smaller than $\alpha$. More about this is said in Sections 4.4 and 4.5.

### *4.4. Spiky priors and false discovery rates*

We commented in Section 4.1 that if the prior density $g(\theta_0)$ is large, it increases the leading term in the expansion for $\delta_n$ (and also $\epsilon_n$) and so it can be expected that spiky priors cause higher false discovery rates. In this section, we address the effect of spiky and flat priors a little more formally.

Consider the general testing problem $H_0 : \theta \leq \theta_0$ vs $H_1 : \theta > \theta_0$, where the natural parameter space $\Omega = (\underline{\theta}, \bar{\theta})$.

Suppose the $\alpha$ $(0 < \alpha < 1)$ level test rejects $H_0$ if $T_n \in C$, where $T_n$ is the test statistic. Let $P_n(\theta) = P_\theta(T_n \in C)$. Let $g(\theta)$ be any fixed density function for $\theta$ and let $g_\tau(\theta) = g(\theta/\tau)/\tau$, $\tau > 0$. Then $g_\tau(\theta)$ is spiky at 0 for small $\tau$ and $g_\tau(\theta)$ is flat for large $\tau$. When $\theta_0 = 0$, under the prior $g_\tau(\theta)$,

$$(34) \qquad \delta_n(\tau) = P(\theta \le 0 | T_n \in C) = \frac{\int_{\underline{\theta}/\tau}^{0} P_n(\tau y) g(y) dy}{\int_{\underline{\theta}/\tau}^{\bar{\theta}/\tau} P_n(\tau y) g(y) dy},$$

and

$$(35) \qquad \epsilon_n(\tau) = P(\theta > 0 | T_n \notin C) = \frac{\int_{\underline{\theta}/\tau}^{0} [1 - P_n(\tau y)] g(y) dy}{\int_{\underline{\theta}/\tau}^{\bar{\theta}/\tau} [1 - P_n(\tau y)] g(y) dy}.$$

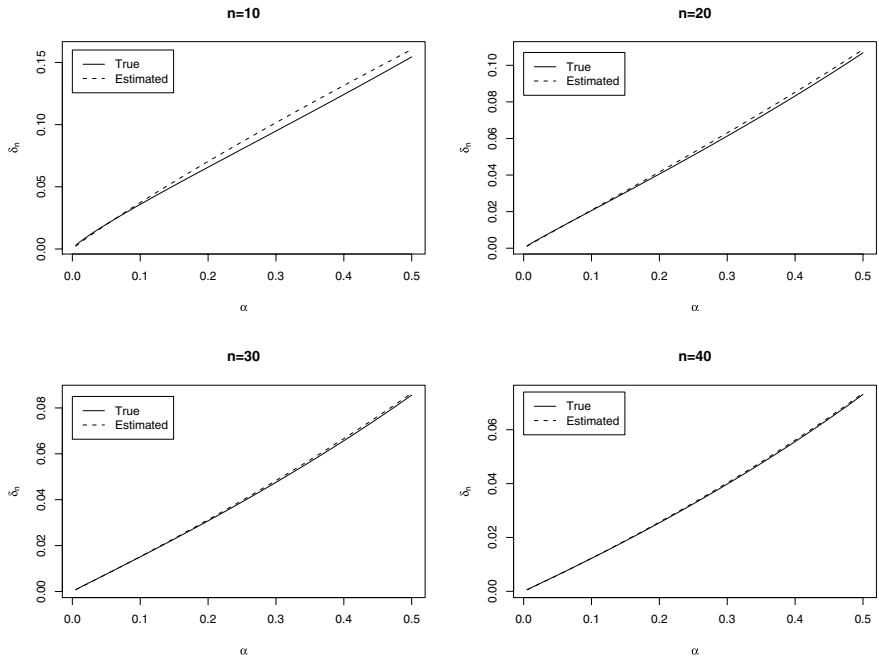

FIG 4. *True and estimated values of $\delta_n$ when the test statistic is the median for the normal-normal case.*

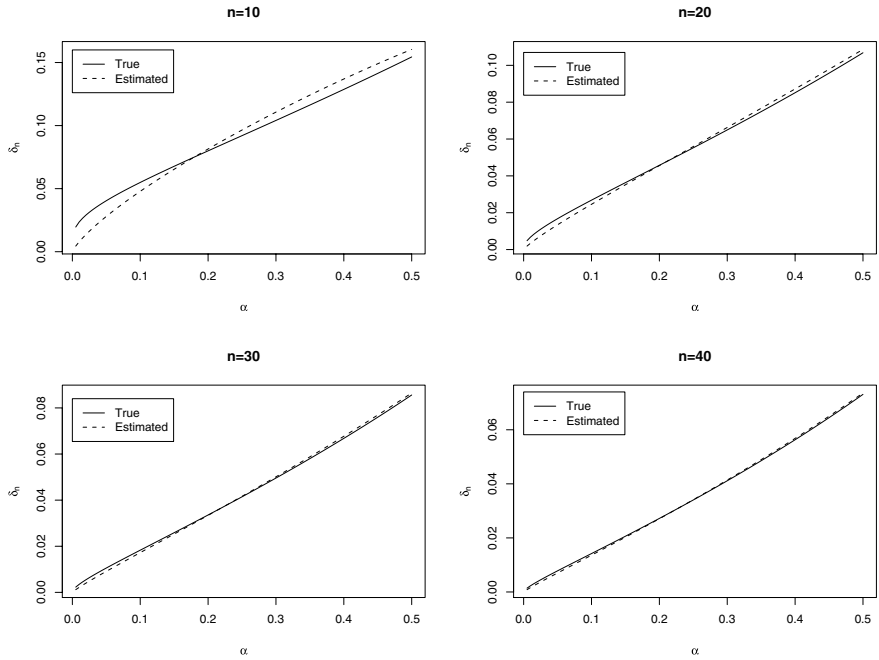

FIG 5. *True and estimated values of $\delta_n$ when the test statistic is the median for the Cauchy-Cauchy case.*

Let as before $\lambda = \int_{-\infty}^{0} g(\theta) d\theta$, the numerator and denominator of (34) be denoted by $A_n(\tau)$ and $B_n(\tau)$ and the numerator and denominator of (35) be denoted by $\tilde{A}_n(\tau)$ and $\tilde{B}_n(\tau)$. Then, we have the following results.

**Proposition 1.** *If* $P_n^-(\theta_0) = \lim_{\theta \to \theta_0-} P_n(\theta)$ *and* $P_n^+(\theta_0) = \lim_{\theta \to \theta_0+} P_n(\theta)$ *both exist and are positive, then*

$$(36) \qquad \lim_{\tau \to 0} \delta_n(\tau) = \frac{\lambda P_n^-(0)}{\lambda P_n^-(0) + (1 - \lambda) P_n^+(0)}$$

*and*

$$(37) \qquad \lim_{\tau \to 0} \epsilon_n(\tau) = \frac{(1 - \lambda)[1 - P_n^+(0)]}{\lambda[1 - P_n^-(0)] + (1 - \lambda)[1 - P_n^+(0)]}.$$

*Proof.* Because $0 \le P_n(\tau y) \le 1$ for all $y$, by simply applying the Lebesgue Dominated Convergence Theorem, $\lim_{\tau \to 0} A_n(\tau) = \lambda P_n^-(0)$, $\lim_{\tau \to 0} B_n(\tau) = \lambda P_n^-(0) + (1-\lambda) P_n^+(0)$, $\lim_{\tau \to 0} \tilde{A}_n(\tau) = (1-\lambda)[1 - P_n^+(0)]$ and $\lim_{\tau \to 0} \tilde{B}_n(\tau) = \lambda[1 - P_n^-(0)] + (1 - \lambda)[1 - P_n^+(0)]$. Substituting in (34) and (35), we get (36) and (37). $\qquad\square$

**Corollary 1.** *If* $0 < \lambda < 1$, $\lim_{\tau \to \infty} P_n(\tau y) = 0$ *for all* $y < 0$, $\lim_{\tau \to \infty} P_n(\tau y) = 1$ *for all* $y > 0$, *then* $\lim_{\tau \to \infty} \delta_n(\tau) = \lim_{\tau \to \infty} \epsilon_n(\tau) = 0$.

*Proof.* Immediate from (36) and (37). $\qquad\square$

It can be seen that $P_n^-(0) = P_n^+(0)$ in most testing problems when the test statistic $T_n$ has a continuous power function. It is true for all the problems we discussed in Sections 3 and 4. If moreover $g(\theta) > 0$ for all $\theta$, then $0 < \lambda < 1$. As a consequence, $\lim_{\tau \to 0} \delta_n(\tau) = \lambda$, $\lim_{\tau \to 0} \epsilon_n(\tau) = 1 - \lambda$, and $\lim_{\tau \to \infty} \delta_n(\tau) = \lim_{\tau \to \infty} \epsilon_n(\tau) = 0$. If $\theta$ is a location parameter, $\theta_0 = 0$ and $g(\theta)$ is symmetric about 0, then $\lim_{\tau \to 0} \delta_n(\tau) = \lim_{\tau \to 0} \epsilon_n(\tau) = 1/2$.

In other words, the false discovery rates are very small for any $n$ for flat priors and roughly 50% for any $n$ for very spiky symmetric priors. This is a qualitatively informative observation.

### 4.5. Pre-experimental promise and post-experimental honesty

We noticed in our example in Section 4.4 that for quite small values of $n$, the post-experimental error rate $\delta_n$ was smaller than the pre-experimental assurance, namely $\alpha$. For any given prior $g$, this is true for all large $n$; but clearly we cannot achieve this uniformly over all $g$, or even large classes of $g$. In order to remain honest, it seems reasonable to demand of a frequentist that $\delta_n$ be smaller than $\alpha$. The question is, typically for what sample sizes can the frequentist assert his honesty.

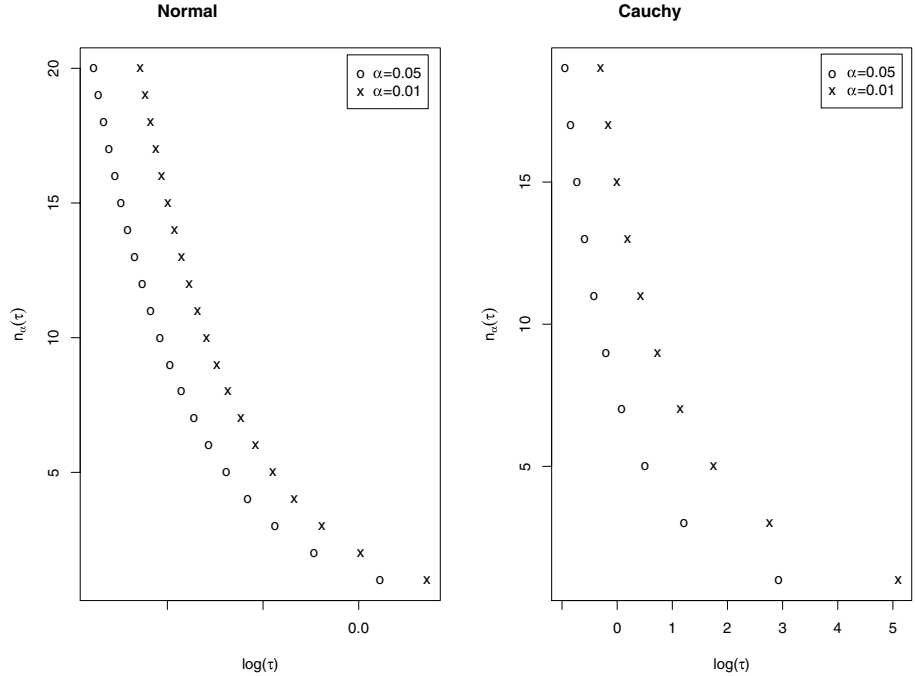

FIG 6. *Plots of $n_\alpha(\tau)$ as functions of $\tau$ for normal-normal test by mean and Cauchy-Cauchy test by median for selected $\alpha$.*

Let us then consider the prior $g_\tau(\theta) = g(\theta/\tau)/\tau$ with fixed $g$, and consider the minimum value of the sample size $n$, denoted by $n_\alpha(\tau)$, such that $\delta_n \leq \alpha$. It can be seen from (36) that $n_\alpha(\tau)$ goes to $\infty$ as $\tau$ goes to 0. This of course was anticipated. What happens when $\tau$ varies from small to large values?

Plots of $n_\alpha(\tau)$ as functions of $\tau$ when the population CDF is $F_\theta(x) = \Phi(x - \theta)$, $g(\theta) = \phi(\theta)$ and the test statistic is $\bar{X}$ are given in the left window of Figure 6. It is seen in the plot that $n_\alpha(\tau)$ is non-increasing in $\tau$ for the selected $\alpha$-values 0.05 and 0.01. Plots of $n_\alpha(\tau)$ when $F_\theta(x) = C(x - \theta)$ and $g(\theta) = c(\theta)$, where $C(\cdot)$ and $c(\cdot)$ are standard Cauchy CDF and PDF respectively, are given in the right window of Figure 6.

In both examples, a modest sample size of $n = 15$ suffices for ensuring $\delta_n \leq \alpha$ if $\tau \geq 1$. For somewhat more spiky priors with $\tau \approx 0.5$, in the Cauchy-Cauchy case, a sample of size just below 30 will be required. In the normal-normal case, even $n = 8$ still suffices.

The general conclusion is that unless the prior is very spiky, a sample of size about 30 ought to ensure that $\delta_n \leq \alpha$ for traditional values of $\alpha$.

## Appendix: Detailed expansions for the exponential family

We now provide the details for the expansions of $A_{n,\theta_{1n}}$ in (13) and $\tilde{A}_{n,\theta_{2n}}$ in (15)

and we also prove that $\underline{R}_{n,\theta_{1n}}$ in (14) and $\bar{R}_{n,\theta_{2n}}$ in (16) are smaller order terms.

Suppose $g(\theta)$ is a three times differentiable proper prior for $\theta$. The expansions are considered for those $\theta_0$ so that the exponential family density has a positive variance at $\theta_0$. Then, we can find two values $\theta_1$ and $\theta_2$ such that $\underline{\theta} < \theta_1 < \theta_0 < \theta_2 < \bar{\theta}$ and the minimum value of $\sigma^2(\theta)$ is positive when $\theta_1 \leq \theta \leq \theta_2$. That is if we let $m_0 = \min_{\theta_1 < \theta < \theta_2} \sigma^2(\theta)$, then $m_0 > 0$. Since $\sigma^2(\theta)$, $k_i(\theta)$, $\rho_i(\theta)$ and $g^{(3)}(\theta)$ are all continuous in $\theta$, each of them is uniformly bounded in absolute value for $\theta \in [\theta_1, \theta_2]$. We denote $M_0$ as the common upper bound of the absolute values of $\sigma^2(\theta)$, $\kappa_i(\theta)$ ($i = 3, 4, 5$), $\rho_i(\theta)$ ($i = 3, 4, 5$), $g(\theta)$, $g'(\theta)$, $g''(\theta)$ and $g^{(3)}(\theta)$.

In the rest of this section, we denote $\theta_{1n} = \theta_0 + (\theta_1 - \theta_0)/n^{1/3}$, $\theta_{2n} = \theta_0 + (\theta_2 - \theta_0)/n^{1/3}$, $x_1 = \sigma_0\sqrt{n}(\theta_1 - \theta_0) - z_\alpha$, $x_2 = \sigma_0\sqrt{n}(\theta_2 - \theta_0) - z_\alpha$, $x_{1n} = \sigma_0\sqrt{n}(\theta_{1n} - \theta_0) - z_\alpha$ and $x_{2n} = \sigma_0\sqrt{n}(\theta_{1n} - \theta_0) - z_\alpha$. As in (13), (14), (15) and (16), we define $A_{n,\theta_{1n}} = P(\theta_{1n} \leq \theta \leq \theta_0, \bar{X} \in C)$, $\underline{R}_{n,\theta_{1n}} = A_n - A_{n,\theta_{1n}}$, $\tilde{A}_{n,\theta_{2n}} = P(\theta_0 < \theta \leq \theta_{2n}, \bar{X} \notin C)$ and $\bar{R}_{n,\theta_{2n}} = \tilde{A}_n - \tilde{A}_{n,\theta_{2n}}$, where $A_n$ and $\tilde{A}_n$ are given by (3) and (4) respectively. We write $B_{n,\theta_1} = P(\theta \geq \theta_{1n}, \bar{X} \in C)$ and $\tilde{B}_{n,\theta_2} = P(\theta \leq \theta_{2n}, \bar{X} \notin C)$. Then, one can also see that $\underline{R}_{n,\theta_{1n}} = B_n - B_{n,\theta_1}$ and $\bar{R}_{n,\theta_{2n}} = \tilde{B}_n - \tilde{B}_{n,\theta_2}$ from definition, where $B_n$ and $\tilde{B}_n$ are given by (5) and (6) respectively.

The following Proposition and Corollary claim that $\underline{R}_{n,\theta_{1n}}$ and $\bar{R}_{n,\theta_{2n}}$ are the smaller order terms. Therefore, the coefficients of the expansions of $A_n$ and $\tilde{A}_n$ are exactly the same as those of $A_{n,\theta_{1n}}$ and $\tilde{A}_{n,\theta_{2n}}$.

**Proposition 2.** *Let $\theta_{1,\tau,n} = \theta_0 + (\theta_1 - \theta_0)/n^\tau$ and $\theta_{2,\tau,n} = \theta_0 + (\theta_2 - \theta_0)/n^\tau$. If $0 \leq \tau < 1/2$, then for any $\ell < \infty$, $\lim_{n\to\infty} n^l \tilde{\beta}_n(\theta_{1,\tau,n}) = \lim_{n\to\infty} n^l [1 - \tilde{\beta}_n(\theta_{2,\tau,n})] = 0$.*

*Proof.* A proof of this can be obtained by simply using Markov's inequality. We omit it. $\qquad\blacksquare$

**Corollary 2.** *For any $l > 0$, $\lim_{n\to\infty} n^l \underline{R}_{n,\theta_{1n}} = \lim_{n\to\infty} n^l \bar{R}_{n,\theta_{2n}} = 0$.*

*Proof.* Since $\tilde{\beta}_n(\theta)$ is nondecreasing in $\theta$, we have

$$n^l \underline{R}_{n,\theta_{1n}} = n^l \int_{\underline{\theta}}^{\theta_{1n}} \tilde{\beta}_n(\theta)g(\theta)d\theta \leq n^l \tilde{\beta}_n(\theta_{1,1/3,n}) \int_{\underline{\theta}}^{\theta_{1n}} g(\theta)d\theta \leq n^l \tilde{\beta}_n(\theta_{1,1/3,n})$$

and similarly $n^l \bar{R}_{n,\theta_{2n}} \leq n^l [1 - \tilde{\beta}_n(\theta_{2,1/3,n})]$. The conclusion is drawn by taking $\tau = 1/3$ in Proposition 2. $\qquad\blacksquare$

In the rest of this section, we will only derive the expansion of $A_{n,\theta_{1n}}$ in detail since the expansion of $\tilde{A}_{n,\theta_{2n}}$ is obtained exactly similarly.

Using the transformation $x = \sigma_0\sqrt{n}(\theta - \theta_0) - z_\alpha$ in the following integral, we have

$$(38) \qquad A_{n,\theta_{1n}} = \frac{1}{\sigma_0\sqrt{n}} \int_{x_{1n}}^{-z_\alpha} \tilde{\beta}_n(\theta_0 + \frac{x + z_\alpha}{\sigma_0\sqrt{n}})g(\theta_0 + \frac{x + z_\alpha}{\sigma_0\sqrt{n}})dx.$$

Note that

$$(39) \qquad \tilde{\beta}_n(\theta_0 + \frac{x + z_\alpha}{\sigma_0\sqrt{n}}) = P_{\theta_0 + \frac{x+z_\alpha}{\sigma_0\sqrt{n}}} \left( \sqrt{n}\frac{\bar{X} - \mu(\theta_0 + \frac{x+z_\alpha}{\sigma_0\sqrt{n}})}{\sigma(\theta_0 + \frac{x+z_\alpha}{\sigma_0\sqrt{n}})} \geq \tilde{k}_{\theta_0,x,n} \right),$$

where

$$(40) \qquad \tilde{k}_{\theta_0,x,n} = \left[ \sqrt{n}\frac{\mu_0 - \mu(\theta_0 + \frac{x+z_\alpha}{\sigma_0\sqrt{n}})}{\sigma_0} + k_{\theta_0,n} \right] \frac{\sigma_0}{\sigma(\theta_0 + \frac{x+z_\alpha}{\sigma_0\sqrt{n}})}.$$

We obtain the coefficients of the expansions of $A_{n,\theta_{1n}}$ in the following steps:

1. The expansion of $g(\theta_0 + \frac{x+z_\alpha}{\sigma_0\sqrt{n}})$ for any fixed $x \in [x_{1n}, -z_\alpha]$ is obtained by using Taylor expansions.
2. The expansion of $\tilde{k}_{\theta_0,x,n}$ for any fixed $x \in [x_{1n}, -z_\alpha]$ is obtained by jointly considering the Cornish-Fisher expansion of $k_{\theta_0,n}$, the Taylor expansion of $\sqrt{n}[\mu_0 - \mu(\theta_0 + \frac{x+z_\alpha}{\sigma_0\sqrt{n}})]/\sigma_0$ and the Taylor expansion of $\sigma_0/\sigma(\theta_0 + \frac{x+z_\alpha}{\sigma_0\sqrt{n}})$.
3. Write the CDF of $\bar{X}$ in the form of $P_\theta[\sqrt{n}\frac{\bar{X}-\mu(\theta)}{\sigma(\theta)} \leq u]$. Formally substitute $\theta = \theta_0 + \frac{x+z_\alpha}{\sigma_0\sqrt{n}}$ and $u = \tilde{k}_{\theta_0,x,n}$ in the Edgeworth expansion of the CDF of $\bar{X}$. An expansion of $\tilde{\beta}_n(\theta_0 + \frac{x+z_\alpha}{\sigma_0\sqrt{n}})$ is obtained by combining it with Taylor expansions for a number of relevant functions (see (47)).
4. The expansion of $A_{n,\theta_{1n}}$ is obtained by considering the product of the expansions of $g(\theta_0 + \frac{x+z_\alpha}{\sigma_0\sqrt{n}})$ and $\tilde{\beta}_n(\theta_0 + \frac{x+z_\alpha}{\sigma_0\sqrt{n}})$ under the integral sign.
5. Finally prove that all the error terms in Steps 1, 2, 3 and 4 are smaller order terms.

We give the expansions in steps 1, 2, 3 and 4 in detail. For the error term study in step 5, we omit the details due to the considerably tedious algebra.

**Step 1**: The expansion of $g(\theta_0 + \frac{x+z_\alpha}{\sqrt{n}})$ is easily obtained by using a Taylor expansion:

$$(41) \quad g(\theta_0 + \frac{x+z_\alpha}{\sigma_0\sqrt{n}}) = g(\theta_0) + g'(\theta_0)\frac{x+z_\alpha}{\sigma_0\sqrt{n}} + \frac{g''(\theta_0)}{2}\frac{(x+z_\alpha)^2}{\sigma_0^2 n} + r_{g,x,n}.$$

where $r_{g,x,n}$ is the error term.

**Step 2:** The Cornish–Fisher expansion of $k_{\theta_0,n}$ ([1], p. 117) is given by

$$(42) \quad k_{\theta_0,n} = z_\alpha + \frac{(z_\alpha^2-1)\rho_{30}}{6\sqrt{n}} + \frac{1}{n}\left[\frac{(z_\alpha^3-3z_\alpha)\rho_{40}}{24} - \frac{(2z_\alpha^3-5z_\alpha)\rho_{30}^2}{36}\right] + r_{1,n},$$

where $r_{1,n}$ is the error term.

The Taylor expansion of the first term inside the bracket of (40) is

$$(43) \quad -(x+z_\alpha) - \frac{\rho_{30}(x+z_\alpha)^2}{2\sqrt{n}} - \frac{\rho_{40}(x+z_\alpha)^3}{6n} + r_{2,x,n}$$

and the Taylor expansion of the term outside of the bracket of (40) is

$$(44) \quad 1 - \frac{\rho_{30}(x+z_\alpha)}{2\sqrt{n}} + \frac{1}{n}\left(\frac{3\rho_{30}^2}{8} - \frac{\rho_{40}}{4}\right)(x+z_\alpha)^2 + r_{3,x,n},$$

where $r_{2,x,n}$ and $r_{3,x,n}$ are error terms.

Plugging (42), (43) and (44) into (40), we get the expansion of $\tilde{k}_{\theta_0,x,n}$ below:

$$(45) \quad \tilde{k}_{\theta_0,x,n} = -x + \frac{1}{\sqrt{n}}f_1(x) + \frac{1}{n}f_2(x) + r_{4,x,n},$$

where $r_{4,x,n}$ is the error term, $f_1(x) = f_{11}x + f_{10}$ and $f_2(x) = f_{23}x^3 + f_{22}x^2 + f_{21}x + f_{20}$, and the coefficients for $f_1(x)$ and $f_2(x)$ are $f_{10} = -(2z_\alpha^2+1)\rho_{30}/6$, $f_{11} = -z_\alpha\rho_{30}/2$, $f_{20} = (z_\alpha^3+2z_\alpha)\rho_{30}^2/9 - (z_\alpha^3+z_\alpha)\rho_{40}/8$, $f_{21} = (7z_\alpha^2/24+1/12)\rho_{30}^2 - z_\alpha^2\rho_{40}/4$, $f_{22} = 0$, $f_{23} = \rho_{40}/12 - \rho_{30}^2/8$.

**Step 3:** The Edgeworth expansion of the CDF of $\bar{X}$ is (Barndorff-Nielsen and Cox ([1], p. 91) and Hall ([11], p. 45)) given below:

$$
\begin{aligned}
&P_{\theta_0 + \frac{x+z_\alpha}{\sigma_0\sqrt{n}}}\left(\sqrt{n}\frac{\bar{X} - \mu(\theta_0 + \frac{x+z_\alpha}{\sigma_0\sqrt{n}})}{\sqrt{\sigma(\theta_0 + \frac{x+z_\alpha}{\sigma_0\sqrt{n}})}} \le u\right) \\
&= \Phi(u) - \frac{\phi(u)}{\sqrt{n}}\frac{(u^2-1)}{6}\rho_3(\theta_0 + \frac{x+z_\alpha}{\sigma_0\sqrt{n}}) - \frac{\phi(u)}{n}[\frac{(u^3-3u)}{24} \\
&\quad \times \rho_4(\theta_0 + \frac{x+z_\alpha}{\sigma_0\sqrt{n}}) + \frac{(u^5-10u^3+15u)}{72}\rho_3^2(\theta_0 + \frac{x+z_\alpha}{\sigma_0\sqrt{n}})] + r_{5,n},
\end{aligned}
\tag{46}
$$

where $r_{5,n}$ is an error term. If we take $\mu = \tilde{k}_{\theta_0,x,n}$ in (46), then the left side is $1 - \tilde{\beta}_n(\theta_0 + \frac{x+z_\alpha}{\sqrt{n}})$ and so

$$
\begin{aligned}
\tilde{\beta}(\theta_0 + \frac{x+z_\alpha}{\sigma_0\sqrt{n}}) &= \Phi(-\tilde{k}_{\theta_0,x,n}) + \frac{\phi(\tilde{k}_{\theta_0,x,n})}{\sqrt{n}}\frac{(\tilde{k}_{\theta_0,x,n}^2-1)}{6}\rho_3(\theta_0 + \frac{x+z_\alpha}{\sigma_0\sqrt{n}}) \\
&\quad + \frac{\phi(\tilde{k}_{\theta_0,x,n})}{n}[\frac{(\tilde{k}_{\theta_0,x,n}^3 - 3\tilde{k}_{\theta_0,x,n})}{24}\rho_4(\theta_0 + \frac{x+z_\alpha}{\sigma_0\sqrt{n}}) \\
&\quad + \frac{(\tilde{k}_{\theta_0,x,n}^5 - 10\tilde{k}_{\theta_0,x,n}^3 + 15\tilde{k}_{\theta_0,x,n})}{72}\rho_3^2(\theta_0 + \frac{x+z_\alpha}{\sigma_0\sqrt{n}})] - r_{5,n}.
\end{aligned}
\tag{47}
$$

Plug the Taylor expansion of $\rho_3(\theta_0 + \frac{x+z_\alpha}{\sigma_0\sqrt{n}})$

$$
\rho_3(\theta_0 + \frac{x+z_\alpha}{\sigma_0\sqrt{n}}) = \rho_{30} + \frac{(x+z_\alpha)}{\sqrt{n}}\left(\rho_{40} - \frac{3}{2}\rho_{30}^2\right) + r_{6,x,n}
\tag{48}
$$

in (47), where $r_{6,x,n}$ is an error term, and then consider the Taylor expansions of the three terms related to $\tilde{k}_{\theta_0,x,n}$ in (47) and also use the expansion (45). On quite a bit of calculations, we obtain the following expansion:

$$
\begin{aligned}
\tilde{\beta}_n(\theta_0 + \frac{x+z_\alpha}{\sigma_0\sqrt{n}}) &= \Phi(x) - \phi(x)\left[\frac{f_1(x)}{\sqrt{n}} + \frac{f_2(x)}{n}\right] - x\phi(x)\left[\frac{f_1(x)}{\sqrt{n}}\right]^2 \\
&\quad + \frac{\phi(x)(x^2-1)}{6\sqrt{n}}[\rho_{30} + \frac{(x+z_\alpha)}{\sqrt{n}}(\rho_{40} - \frac{3}{2}\rho_{30}^2)] + \frac{\rho_{30}}{6n}\phi(x)(x^3-3x)f_1(x) \\
&\quad + \frac{\phi(x)}{n}[\frac{(x^3-3x)}{24}\rho_{40} + \frac{(x^5-10x^3+15x)}{72}\rho_{30}^2] + r_{7,x,n} \\
&= \Phi(x) + \frac{\phi(x)}{\sqrt{n}}g_1(x) + \frac{\phi(x)}{n}g_2(x) + r_{7,x,n},
\end{aligned}
\tag{49}
$$

where $r_{7,x,n}$ is an error term, $g_1(x) = g_{12}x^2 + g_{11}x + g_{10}$, $g_2(x) = g_{20} + g_{21}x + g_{22}x^2 + g_{23}x^3 + g_{24}x^4 + g_{25}x^5$, and the coefficients of $g_1(x)$ and $g_2(x)$ are $g_{12} = \rho_{30}/6$, $g_{11} = z_\alpha\rho_{30}/2$, $g_{10} = z_\alpha^2\rho_{30}/3$, $g_{25} = \rho_{30}^2/72$, $g_{24} = -z_\alpha\rho_{30}^2/12$, $g_{23} = \rho_{40}/8 - 13z_\alpha^2\rho_{30}^2/72 - 7\rho_{30}^2/24$, $g_{22} = z_\alpha\rho_{40}/6 - z_\alpha^3\rho_{30}^2/6 - z_\alpha\rho_{30}^2/12$, $g_{21} = (z_\alpha^2/4 - 7/24)\rho_{40} - z_\alpha^4\rho_{30}^2/18 - 13z_\alpha^2\rho_{30}^2/72 + 4\rho_{30}^2/9$, $g_{20} = (z_\alpha^3/8 - z_\alpha/24)\rho_{40} - (z_\alpha^3/9 - z_\alpha/36)\rho_{30}^2$.

**Step 4:** The expansion of $A_{n,\theta_{1n}}$ is obtained by plugging the expansions of $\tilde{\beta}(\theta_0 + \frac{x+z_\alpha}{\sigma_0\sqrt{n}})$ and $g(\theta_0 + \frac{x+z_\alpha}{\sigma_0\sqrt{n}})$. On careful calculations,

$$
A_{n,\theta_{1n}} = \frac{a_1}{\sqrt{n}} + \frac{a_2}{n} + \frac{a_3}{n^{3/2}} + r_{8,n},
\tag{50}
$$

where $r_{8,n}$ is an error term, $a_1 = (g(\theta_0)/\sigma_0)[\phi(z_\alpha) - \alpha z_\alpha]$, $a_2 = \rho_{30} g(\theta_0)[\alpha + 2\alpha z_\alpha^2 - 2z_\alpha \phi(z_\alpha)]/(6\sigma_0) - g'(\theta_0)[\alpha(z_\alpha^2 + 1) - z_\alpha \phi(z_\alpha)]/(2\sigma_0^2)$, and $a_3 = [h_{11}\phi(z_\alpha) + \alpha h_{12}][(g''(\theta_0)/(6\sigma_0^3)] + [h_{21}\phi(z_\alpha) + \alpha h_{22}][g'(\theta_0)/\sigma_0^2] + [h_{31}\phi(z_\alpha) + \alpha h_{32}][g(\theta_0)/\sigma_0]$, where $h_{11} = z_\alpha^2 + 2$, $h_{12} = -(z_\alpha^3 + 3z_\alpha)$, $h_{21} = -(\rho_{30}/3)(z_\alpha^2 + 1)$, $h_{22} = (\rho_{30}/3)(z_\alpha^3 + 2z_\alpha)$, $h_{31} = -z_\alpha^4 \rho_{30}^2/36 + 4z_\alpha^2 \rho_{30}^2/9 + \rho_{30}^2/36 - 5z_\alpha^2 \rho_{40}/24 + \rho_{40}/24$, $h_{32} = -5z_\alpha^3 \rho_{30}^2/18 - 11z_\alpha \rho_{30}^2/36 + z_\alpha^3 \rho_{40}/8 + z_\alpha \rho_{40}/8$. These $a_1$, $a_2$ and $a_3$ are the coefficients in the expansion of (23).

The computation of the coefficients of the expansions of $A_{n,\theta_{1n}}$ is now complete. The rest of the work is to prove that all the error terms are smaller order terms. But first we give the results for the expansion of $\tilde{A}_{n,\theta_{2n}}$. The details for the expansions of $\tilde{A}_{n,\theta_{2n}}$ are omitted.

**Expansion of $\tilde{A}_{n,\theta_{2n}}$:** The expansion of $\tilde{A}_{n,\theta_{2n}}$ can be obtained similarly by simply repeating all the steps for $A_{n,\theta_{1n}}$. The results are given below:

$$(51) \qquad \tilde{A}_{n,\theta_{2n}} = \frac{\tilde{a}_1}{\sqrt{n}} + \frac{\tilde{a}_2}{n} + \frac{\tilde{a}_3}{n^{3/2}} + r_{9,n},$$

where $r_{9,n}$ is an error term, $\tilde{a}_1 = g(\theta_0)[\phi(z_\alpha) + (1 - \alpha)z_\alpha]/\sigma_0$, $\tilde{a}_2 = g'(\theta_0)[(1 - \alpha)(z_\alpha^2 + 1) + z_\alpha \phi(z_\alpha)]/(2\sigma_0^2) - \rho_{30} g(\theta_0)[(1 - \alpha)/6 + (1 - \alpha)z_\alpha^2/3 + z_\alpha \phi(z_\alpha)/3]/\sigma_0$, and $\tilde{a}_3 = (g''(\theta_0)/6\sigma_0^3)[h_{11}\phi(z_\alpha) - (1 - \alpha)h_{12}] - (g'(\theta_0)/\sigma_0^2)[-h_{21}\phi(z_\alpha) + (1 - \alpha)h_{22}] - (g(\theta_0)/\sigma_0)[-h_{31}\phi(z_\alpha) + (1 - \alpha)h_{32}]$, where $h_{11}$, $h_{12}$, $h_{21}$, $h_{22}$, $h_{31}$ and $h_{32}$ are the same as defined in Step 3. These $\tilde{a}_1$, $\tilde{a}_2$ and $\tilde{a}_3$ are the coefficients in the expansion of (24).

**Remark.** The coefficients of expansions of $\delta_n$ and $\epsilon_n$ are obtained by simply using formula (7) with $a_1$, $a_2$ and $a_3$ in (23) and also the coefficient $\tilde{a}_1$, $\tilde{a}_2$ and $\tilde{a}_3$ in (24) respectively.

**Step 5:** (Error term study in the expansions of $A_{n,\theta_{1n}}$). We only give the main steps because the details are too long. Recall from equation (38) that the range of integration corresponding to $A_{n,\theta_{1n}}$ is $x_{1n} \le x \le -z_\alpha$. In this case, we have $\lim_{n\to\infty} x_{1n} = \infty$ and $\lim_{n\to\infty} x_{1n}/\sqrt{n} = -z_\alpha$. This fact is used when we prove the error term is still a smaller order term when we move it out of the integral sign.

(I)   In (41), since $g^{(3)}(\theta)$ is uniformly bounded in absolute values, $r_{g,x,n}$ is absolutely bounded by a constant times $n^{-3/2}(x + z_\alpha)^2$

(II)  From Barndorff-Nielsen and Cox [4.5, pp 117], the error term $r_{1,n}$ in (42) is absolutely uniformly bounded by a constant times $n^{-3/2}$.

(III) In (43) and (44), since $\rho_i(\theta)$ and $\kappa_i(\theta)$ ($i = 3, 4, 5$) are uniformly bounded in absolute values, the error term $r_{2,x,n}$ is absolutely bounded by a constant times $n^{-3/2}(x + z_\alpha)^4$ and the error term $r_{3,x,n}$ is absolutely bounded by a constant times $n^{-3/2}(x + z_\alpha)^3$.

(IV)  The exact form of the error term $r_{4,x,n}$ in (45) can be derived by considering the higher order terms and their products in (42), (43) and (44) for the derivation of expression (45). The computation is complicated but straightforward. However, still, since $\rho_i(\theta)$ and $\kappa_i(\theta)$ ($i = 3, 4, 5$) are uniformly bounded in absolute values, $r_{4,x,n}$ is absolutely bounded by $n^{-3/2}P_1(|x|)$, where $P_1(|x|)$ is a seventh degree polynomial and its coefficients do not depend on $n$.

(V)   Again, from Barndorff-Nielsen and Cox ([1], p. 91), the error term $r_{5,n}$ in (46) is absolutely bounded by a constant times $n^{-3/2}$.

(VI)  The error term $r_{6,x,n}$ in (48) is absolutely bounded by a constant times $n^{-1}(x + z_\alpha)^2$ since $\rho_i(\theta)$ and $\kappa_i(\theta)$ ($i = 3, 4, 5$) are uniformly bounded in absolute values.

(VII)  This is the critical step for the error term study since we need to prove that the error term is still a smaller order term when it is moved out of the integral in (50). We need to study the behaviors of $\Phi(-\tilde{k}_{\theta_0,x,n})$ and $\phi(\tilde{k}_{\theta_0,x,n})$ as $n \to \infty$ for all $x \in [x_{1n}, -z_\alpha]$ uniformly (see (49) in detail). This also explains why we choose $\theta_{1n} = \theta_0 + (\theta_1 - \theta_0)/n^{1/3}$ and $x_{1n} = \sigma_0\sqrt{n}(\theta_{1n} - \theta_0) - z_\alpha$ at the beginning of this section, since in this case $|\tilde{k}_{\theta_0,x,n} + x|$ is uniformly bounded by $|x|/2 + 1$ for a sufficiently large $n$. Then for sufficiently large $n$, the error term $|r_{7,x,n}|$ in (49) is uniformly bounded by $|r_{7,x,n}| \leq \phi(x/2 + 1)P_2(|x|)$ where $P_2(|x|)$ is a twelveth degree polynomial of $|x|$ and its coefficients do not depend on $n$.

(VIII)  Finally, we can show that the error term $r_{8,n}$ in (50) in $O(n^{-2})$. This is tedious but straightforward. It is proven by considering each of the ten terms in $r_{8,n}$ separately.

**Remark.** We can similarly prove that the error term $r_{9,n}$ in (51) corresponding to $\tilde{A}_{n,\theta_{2n}}$ is $O(n^{-2})$. Since the steps are very similar, we do not mention them.

## Acknowledgment

It is a pleasure to thank Michael Woodroofe for the numerous conversations we had with him on the foundational questions and the technical aspects of this paper. This paper would never have been written without Michael's association. We would also like to thank two anonymous referees for their thoughtful remarks and Jiayang Sun for her editorial input. We thank Larry Brown for reading an earlier draft of the paper and for giving us important feedback and to Persi Diaconis for asking a question that helped us in interpreting the results.